\documentclass[12pt]{article}

\usepackage{amsfonts}
\usepackage{latexsym}

\def\medskipamount{8pt} 
\def\arraystretch{1.3}
\begin{document}

% Changing the way equations are numbered

\renewcommand{\thesection}%
{Part \Roman{section}\hspace{.2ex}:\hspace{-1ex}}
\renewcommand{\thesubsection}{\arabic{subb}}

\newcounter{subb}
\setcounter{subb}{0}

\newcommand{\bit}[1]{\pagebreak[3] \addtocounter{subb}{1}
  \subsection{\hspace{-2ex} #1}\setcounter{equation}{0}\penalty12000}
\renewcommand{\theequation}{\thesubsection.\arabic{equation}}

% Putting references in parentheses and BOLDFACE

\newcommand{\re}[1]{\mbox{\bf (\ref{#1})}}

% Putting equation numbers in boldface and on the left

\catcode`\@=\active
\catcode`\@=11
\def\@eqnnum{\hbox to .01pt{}\rlap{\bf \hskip -\displaywidth(\theequation)}}
\catcode`\@=12

%%%%%%%%%%%%%%%%%%%%%%%%%%%%%%%%%%%%%%%%%%%%%%%%%%%%%%%%%%%%%%%%%%%

% NEW COMMANDS

\catcode`\@=\active
\catcode`\@=11
\newcommand{\nc}{\newcommand}

% "Environments" with boldface titles

   % To begin a slanted environment
\nc{\bs}[1]{\addpenalty{-1000}
\addvspace{\medskipamount}\penalty10000
\refstepcounter{equation}
\noindent \begin{em}{\bf (\theequation)} \nopagebreak}

   % To begin a slanted environment with a list
\nc{\bsitem}[1]{\em\refstepcounter{equation}
\hskip -\leftmargin\rlap {\bf  (\theequation)}\hskip\leftmargin}

   % To end a slanted environment
\nc{\es}{\par \end{em} \addvspace{\medskipamount} }

   % To begin a roman environment
\nc{\br}[1]{\addvspace{\medskipamount} %\pagebreak[3]
\refstepcounter{equation}
\noindent {\bf (\theequation) #1.} \nopagebreak }

   % To end a roman environment
\nc{\er}{\par \addvspace{\medskipamount} }

   % To begin a list with roman numeral indices
\newcounter{index}
\nc{\bl}{\begin{list}{{\rm (\roman{index})}}{\usecounter{index}}}

   % To end such a list
\nc{\el}{\end{list}}

   % To begin a proof
\nc{\pf}{\addvspace{\medskipamount} \par \noindent {\em Proof}}

   % To end a proof
\nc{\fp}{\phantom{.} \hfill \mbox{     $\Box$} \par
\addvspace{\medskipamount}}

\nc{\beq}{\begin{equation}}
\nc{\eeq}{\end{equation}}

\nc{\beqas}{\begin{eqnarray*}}
\nc{\eeqas}{\end{eqnarray*}}

% Blackboard boldface

\nc{\C}{\mathbb C}
\nc{\Hyp}{\mathbf H}
\nc{\Pj}{\mathbb P}
\nc{\Q}{\mathbb Q}
\nc{\R}{\mathbf R}
\nc{\Z}{\mathbb Z}

% Roman operators

\nc{\oper}[1]{\mathop{\mathchoice{\mbox{\rm #1}}{\mbox{\rm #1}}
{\mbox{\rm \scriptsize #1}}{\mbox{\rm \tiny #1}}}\nolimits}
\nc{\Ad}{\oper{Ad}}
\nc{\Aut}{\oper{Aut}}
\nc{\const}{\oper{const.}}
\nc{\ch}{\oper{ch}}
\nc{\diag}{\oper{diag}}
\nc{\Diff}{\oper{Diff}}
\nc{\End}{\oper{End}}
\nc{\END}{\oper{$\mathbf{End}$}}
\nc{\Gr}{\oper{Gr}}
\nc{\Hom}{\oper{Hom}}
\nc{\HOM}{\oper{$\mathbf{Hom}$}}
\nc{\id}{\oper{id}}
\nc{\Jac}{\oper{Jac}}
\nc{\Map}{\oper{Map}}
\nc{\mod}{\oper{mod}}
\nc{\Pol}{\oper{Pol}}
\nc{\Quot}{\oper{Quot}}
\nc{\slope}{\oper{slope}}
\nc{\Sym}{\oper{Sym}}
\nc{\rk}{\oper{rk}}
\nc{\td}{\oper{td}}
\nc{\tr}{\oper{tr}}

\nc{\operlimits}[1]{\mathop{\mathchoice{\mbox{\rm #1}}{\mbox{\rm #1}}
{\mbox{\rm \scriptsize #1}}{\mbox{\rm \tiny #1}}}}

\nc{\Coeff}{\operlimits{Coeff}}
\nc{\Res}{\operlimits{Res}}

\nc{\medoplus}{\operlimits{$\bigoplus$}}

% Greek abbreviations

\nc{\al}{\alpha}
\nc{\be}{\beta}
\nc{\ep}{\varepsilon}
\nc{\ga}{\gamma}
\nc{\Ga}{\Gamma}
\nc{\la}{\lambda}
\nc{\La}{\Lambda}
\nc{\si}{\sigma}
\nc{\Sig}{{\Sigma}}
\nc{\Om}{\Omega}

% Miscellaneous

\nc{\elll}{l}
\nc{\well}{w}
\nc{\tag}{}

% Downward arrows with labels
\nc{\down}[1]{{\phantom{\scriptstyle #1}
\hbox{$\left\downarrow\vbox to
    9.5pt{}\right.\nulldelimiterspace=0pt \mathsurround=0pt$}
\raisebox{.4ex}{$\scriptstyle #1$}}}
\nc{\leftdown}[1]{
\hbox{$\raisebox{.4ex}{$\scriptstyle #1$}\left\downarrow\vbox to
    9.5pt{}\right.\nulldelimiterspace=0pt \mathsurround=0pt$}
    \phantom{\scriptstyle #1}}
\nc{\bigdowneq}{{\Big| \! \Big|}}

% Calligraphic letters
\nc{\ca}{{\mathcal A}}
\nc{\ce}{{\mathcal E}}
\nc{\ci}{{I}}
\nc{\cl}{{\mathcal L}}
\nc{\co}{{\mathcal O}}
\nc{\cu}{{\mathcal U}}
\nc{\cv}{{\mathcal V}}
\nc{\B}{{\mathcal B}}
\nc{\F}{{\mathcal F}}
\nc{\G}{{\mathcal G}}
\renewcommand{\H}{{\mathcal H}}
\nc{\M}{{\mathcal M}}
\nc{\N}{{\mathcal N}}

% Lie groups
\nc{\GL}[1]{{{\rm GL(}#1{\rm )}}}
\nc{\PGL}[1]{{{\rm PGL(}#1{\rm )}}}
\nc{\SL}[1]{{{\rm SL(}#1{\rm )}}}
\nc{\SU}[1]{{{\rm SU(}#1{\rm )}}}
\nc{\U}[1]{{{\rm U(}#1{\rm )}}}
\nc{\Sp}[1]{{{\rm Sp(}#1{\rm )}}}
\nc{\cx}{\C^\times}

% Reasonable sized binomial coefficients
\nc{\bino}[2]{\mbox{\Large $#1 \choose #2$}}

% Fractions
\nc{\half}{\mathchoice{{\textstyle \frac{\scriptstyle 1}{\scriptstyle
2}}} {{\textstyle\frac{\scriptstyle 1}{\scriptstyle 2}}}
{\frac{\scriptscriptstyle 1}{\scriptscriptstyle 2}}
{\frac{\scriptscriptstyle 1}{\scriptscriptstyle 2}}}
\nc{\quarter}{\mathchoice{{\textstyle \frac{\scriptstyle
1}{\scriptstyle 4}}} {{\textstyle\frac{\scriptstyle 1}{\scriptstyle
4}}} {\frac{\scriptscriptstyle 1}{\scriptscriptstyle 4}}
{\frac{\scriptscriptstyle 1}{\scriptscriptstyle 4}}}
\nc{\ratio}[2]{\mathchoice{ {\textstyle \frac{\scriptstyle
#1}{\scriptstyle #2}}} {{\textstyle\frac{\scriptstyle #1}{\scriptstyle
#2}}} {\frac{\scriptscriptstyle #1}{\scriptscriptstyle #2}}
{\frac{\scriptscriptstyle #1}{\scriptscriptstyle #2}}}

% To stack 2 lines in a limit of summation
% (Use \scriptsize inside \mbox if non-math test is needed)
\nc{\stack}[2]{
{\def\arraystretch{.7}\def\arraycolsep{0pt}\begin{array}{c}
\scriptstyle #1 \\ \scriptstyle #2 \end{array}\def\arraystretch{1}} }

% Reasonable sized parentheses to counteract bug in 12pt LaTeX

\nc{\Left}[1]{\hbox{$\left#1\vbox to
    11.5pt{}\right.\nulldelimiterspace=0pt \mathsurround=0pt$}}
\nc{\Right}[1]{\hbox{$\left.\vbox to
    11.5pt{}\right#1\nulldelimiterspace=0pt \mathsurround=0pt$}}
\nc{\LEFT}[1]{\hbox{$\left#1\vbox to
    15.5pt{}\right.\nulldelimiterspace=0pt \mathsurround=0pt$}}
\nc{\RIGHT}[1]{\hbox{$\left.\vbox to
    15.5pt{}\right#1\nulldelimiterspace=0pt \mathsurround=0pt$}}

\nc{\lp}{\raisebox{-.1ex}{\rm\large(}\hspace{-.2ex}}
\nc{\rp}{\hspace{-.2ex}\raisebox{-.1ex}{\rm\large)}}

% Truly miscellaneous
\nc{\emb}{\hookrightarrow}
\nc{\jay}{j}
\nc{\lrow}{\longrightarrow}
\nc{\llrow}{\hbox to 25pt{\rightarrowfill}}
\nc{\lllrow}{\hbox to 29pt{\rightarrowfill}}
\nc{\mbar}{\overline{M}}
\nc{\dbar}{\overline{\partial}}
\nc{\sans}{\,\backslash\,}
\nc{\st}{\, | \,}

\nc{\bee}{{\bf E}}
\nc{\bphi}{{\bf \Phi}}

\nc{\rs}{C}

\nc{\T}{T}
\nc{\Tor}{{\mathcal T}}

\nc{\g}{{\mathfrak g}}

\catcode`\@=12

%%%%%%%%%%%%%%%%%%%%%%%%%%%%%%%%%%%%%%%%%%%%%%%%%%%%%%%%%%%%%%%%

% Left-justified title

\noindent {\LARGE \bf Relations in the cohomology ring
\smallskip
\\of the moduli space of rank 2 Higgs bundles}\smallskip\\
%{\Large \bf II: Relations in rank 2}
\smallskip \\
{\bf Tam\'as Hausel }\\
Department of Mathematics, University of California,
Berkeley, Calif. 94720\smallskip \\
{\bf Michael Thaddeus } \\
Department of Mathematics, Columbia University,
New York, N.Y. 10027
\renewcommand{\thefootnote}{}
\footnotetext{T.H. supported by NSF grant DMS--97--29992;
  M.T. supported by NSF grant DMS--98--08529.}

\bigskip

\noindent Let $\rs$ be a smooth complex projective curve of genus
$g$. By a {\em Higgs bundle} on $\rs$ we shall mean a pair
$(E,\phi)$ consisting of a holomorphic vector bundle $E$ on $\rs$
and a section $\phi \in H^0\left(C,\End E \otimes K(np)\right)$,
where $n \geq 0$ is fixed and $p \in \rs$ is a basepoint.  With
the appropriate notion of stability, there exists a
quasi-projective moduli space $\H_n$ of stable Higgs bundles of
fixed rank and degree. The aim of this paper is to characterize
the rational cohomology ring of $\H_n$ when the rank is 2 and the
degree is odd.

In fact, we have given a complete set of generators for this ring
in another paper \cite{ht}.  So it is now a question of
finding the relations between these generators.  Even though this
is a natural companion to the generation problem, the ideas and
methods with which it is studied have quite a different flavor. In
particular, there is much more explicit calculation.

What makes the rank 2 case tractable is that the number of
generators is manageable: just $\ep_1, \dots, \ep_{2g} \in H^1$,
$\al \in H^2$, $\psi_1, \dots, \psi_{2g} \in H^3$, and $\be \in
H^4$.  In arbitrary rank, the number of generators is quite large,
and finding all the relations seems out of reach at the moment. At
any rate, it has not even been done for the compact moduli space
$\N$ of stable bundles which lies inside $\H_n$.

Now in rank 2 the relations in the cohomology ring of $\N$ have been
calculated by several authors \cite{bar,kn,st,zag}, and of course our
ideal of relations must contain theirs.  However, the answer to our
problem is actually more explicit in that the relations are given by a
simple closed formula, rather than by a recursion in the genus $g$.

Yet in another sense our answer is more complicated.  Although
our relations are completely explicit, there are very many of them,
and there is no escaping this.  The so-called ``invariant part'' of
$H^*(\N)$, for example, requires only 3 relations to generate all the
others multiplicatively, while the corresponding number for
$H^*(\H_n)$ grows quadratically in $g$.

Within this large crowd, some old friends stand out.  For example,
$\be^g$, conjectured by Newstead \cite{new} to be a relation on $\N$
in 1972 and shown to be one much later \cite{kir2,cft}, remains a
relation on $\H_0$, though not on all of the $\H_n$.

Another familiar part of the story is a proposition proved in
\S\ref{divisible}, stating that, if $\rho$ is a relation on $\H_0$
at genus $g-1$, then $\psi_j\psi_{j+g} \rho$ is a relation on
$\H_0$ at genus $g$.  Results like this have appeared several
times before in the subject \cite{munoz,cft,thad2}.  They rely
strongly on the interpretation of the moduli space in terms of
flat connections (or more precisely, connections of constant
central curvature).  Indeed, this marks the sole point where this
interpretation, otherwise banished in favor of Higgs bundles, is
briefly recalled.

The contents of the paper may be summarized as follows.  Section
\ref{prelim} recalls the existing results we need, notably the
generation theorem.  Section \ref{state} states the main result of
the paper, which gives an explicit basis $\rho^c_{r,s,t}$ for the
ideal of relations in rank 2.  Section \ref{equivariant} then reviews
some basic facts about equivariant cohomology, which is the main
technical tool.  Section \ref{fixed} classifies the components of the
fixed-point set of the $\cx$-action on $\H_n$ in whose equivariant
cohomology we will work.  These are closely related to the symmetric
products of $\rs$, and so \S\ref{symmetric} reviews some relevant
facts about symmetric products of curves.  Using these facts,
\S\ref{restriction} computes the restrictions of the generators to
each component of the fixed-point set.  This is crucial information if
a polynomial in these classes is to be shown a relation, since a
theorem of Kirwan \re{x}(iii) implies that a polynomial in the
generators is an equivariant relation if and only if its restriction
to every component of the fixed-point set vanishes.

Section \ref{polynomials} defines some polynomials $\xi^k_{r,s}$ in
$\al, \be, \ga$, by a recursive process.  Then in \S\ref{relations} we
finally compute our first equivariant relations, by showing that
certain combinations of the $\xi^k_{r,s}$ vanish when restricted to
the fixed-point set.  These relations are expressed in terms of the
explicit polynomials $\rho^c_{r,s,t}$ by a purely algebraic argument,
given in \S\ref{algebra}.  The proposition discussed before, about
multiplication by $\ga$, is proved next, in \S\ref{divisible}; it
provides many new relations which are divisible by $\ga$.

Up to this point, all the relations we have considered are polynomials
in $\al, \be, \ga$, and thus in particular invariant under the action
of the symplectic group $\Ga = \Sp{2g,\Z}$; \S\ref{gamma} remedies
this situation, showing how the invariant relations may be used to
deduce many more non-invariant relations.  In fact, the relations on
$\H_n$ at genus $g$ divisible by $\psi_1\psi_2 \cdots \psi_k$ turn out
to be
exactly the invariant relations on $\H_{n+k}$ at genus $g-k$.
Thus, even if one is interested only in the space $\H_0$ parametrizing
connections of constant central curvature, one needs to understand the
cohomology of $\H_n$ for $n > 0$.

Finally, \S\ref{wrap} wraps up the proof of the main theorem, by
showing that the number of relations we already have in hand equals
the total number.

The story of how these relations were conjectured and proved is worth
telling briefly.  It was clear from the outset that Kirwan's theorem
would be an invaluable tool; that it could be used, as in
\S\ref{relations}, to decide whether a polynomial in the equivariant
generators is a relation; and that this was, in principle, completely
algorithmic.  The first author was able to implement this algorithm in
the computer software package Macaulay 2 \cite{gs} and crank out the
equivariant relations for $g \leq 7$.  Their restrictions to ordinary
cohomology showed some discernible patterns, so after two weeks in
Oberwolfach, we were able to guess all of the ordinary relations for
general $g$, in roughly the form $\rho^c_{r,s,t}$ given in \re{mmm}.  It
remained only to guess their equivariant extensions, that is, the
relations in the equivariant cohomology, which has one extra
generator.  Guessing these equivariant extensions took the better part
of two years.

More precisely, we were {\em never} able to guess a closed formula or
generating function for the equivariant extensions of the
$\rho^c_{r,s,t}$ themselves, even though we knew them in genus $\leq
7$.  If we could have, this would have been a much shorter paper.  We
hope that someone will guess these extensions some day; what we
actually did, though, was to find rather different equivariant
relations \re{vv}, and then show by several artful maneuvers,
including the proposition mentioned above, that these imply the
relations $\rho^c_{r,s,t}$ we originally found.

Some other papers of the authors have explored different aspects of
the theory of Higgs bundles: for example, the compactification of the
moduli space \cite{hausel1}, its intersection numbers in the compactly
supported cohomology \cite{hausel2}, and the upward and downward flows
from the components of the fixed-point set \cite{thad3}.  Some
relationships between these topics and the contents of this paper are
sketched in the last section, \S\ref{tie-up}.

{\em Notation and conventions.}  Throughout the paper, $\rs$
denotes the smooth projective curve of genus $g$ over which we work.
Its cohomology has the usual generators $e_1, \dots, e_{2g} \in H^1$,
and $\si = e_j e_{j+g} \in H^2$.  The $m$th symmetric product of $\rs$
is denoted $\rs_m$, and the Jacobian of degree $d$ line bundles on
$\rs$ is denoted $\Jac^d \rs$.  The letters $\N$, $\H_n$, and $\M_n$
denote moduli spaces over $\rs$, respectively, of stable bundles $E$
having $\La^2 E$ isomorphic to a fixed line bundle $\Xi$, of Higgs
bundles with values in $K(n) = K \otimes \co(np)$, and Higgs bundles
$(E, \phi)$ with values in $K(n)$ having $\La^2 E \cong \Xi$ and $\tr
\phi = 0$.  Groups are denoted $T = \cx$, $\Ga = \Sp{2g, \Z}$, and
$\Sig = \Z_2^{2g}$.

We use the term {\em total degree} to mean {\em half\/} the
ordinary degree of a cohomology class.  The notation $(\phantom{.})_d$
means the part of a cohomology class in total degree $d$.  All
cohomology is with rational coefficients unless otherwise stated.

We do {\em not} assume $g \geq 2$: the moduli spaces $\M_0$, $\H_0$
and $\N$ are trivial or empty if $g=0$ or $1$, but $\M_n$ and $\H_n$
for $n>0$ are not so trivial, and they play an important role even for
understanding $g \geq 2$.

{\em Acknowledgements.}  We are grateful to the Mathe\-mat\-isches
Forschungs\-institut Ober\-wolfach for its kind hospitality during
three weeks in the summer of 1998, when the ideas for the present
paper took shape.  Computer calculations performed there in the
software package Macaulay 2, by Daniel Grayson and Michael Stillman,
were indispensable in formulating the main result.  So were subsequent
calculations in both Macaulay and Maple.  We are also very
grateful to Andrew Kresch for directing us to the work of Shalosh B.
Ekhad, and to Ekhad for pointing out the recurrence relations that
enabled us to solve the binomial coefficient identities of
\S\ref{algebra}.  We thank Don Zagier for kindly sending us the proof
of a delightful generating function formula from his paper \cite{zag},
which became \re{ss}.  We both wish to thank Nigel Hitchin for his
advice and encouragement, and Tara Brendle and Peter Newstead for
helpful discussions.  Finally, the first author thanks the Institute
for Advanced Study and the Miller Institute for Basic Research in
Science for their hospitality and support during the academic years
1998--99 and 1999--2000, respectively.

\bit{Preliminaries and review}
\label{prelim}

Let $\rs$ be a smooth complex projective curve of genus $g$, and let
$p \in \rs$ be a distinguished point.  For $n \geq 0$, denote by
$K(n)$ the line bundle $K \otimes \co(np)$, where $K$ is the canonical
bundle.  A {\em Higgs bundle} with values in $K(n)$ is a pair $(E,
\phi)$ consisting of a holomorphic vector bundle $E$ over $\rs$ and a
{\em Higgs field} $\phi \in H^0(\End E \otimes K(n))$.  It is {\em
stable} if for all proper subbundles $F \subset E$ such that
$\phi(F) \subset F \otimes K(n)$, $\deg F / \rk F < \deg E / \rk E$.

The existing results we shall need to recall are few and, except
for the authors' result stated as \re{four} below, have 
been known for some time.  They can be summarized as follows.

First, there are some elementary facts about stable Higgs bundles
\cite[4.2,4.3]{ht}. 

\bs{Theorem}
\label{one}
If $X$ parametrizes two families $(\bee, \bphi)$ and $(\bee', \bphi')$
of stable Higgs bundles, and
$(\bee^{\phantom{.}}_x,\bphi^{\phantom{.}}_x) \cong (\bee_x',\bphi_x')$
for all $x\in X$, then $\bee' = \bee \otimes L$ for some line bundle
$L$ over $X$, and $\bphi' = \bphi$.  Moreover, if an action of $\cx$
on $X$ lifts to $\bee$ and $\bee'$ preserving the Higgs fields, then
it lifts to $L$ also so that $\bee' = \bee \otimes L$ equivariantly.
\es

Next, there is a moduli space, constructed by Simpson \cite{simp} and
Nitsure \cite{nit}.

\bs{Theorem}
\label{two}
There exists a moduli space $\H_n$ of stable Higgs bundles of rank 2
and degree 1 with values in $K(n)$, which is a smooth quasi-projective
variety.  It admits a universal family $(\bee, \bphi)$, and the
$\cx$-action on $\H_n$ given by $\la(E\phi) = (E, \la\phi)$ lifts to
this family.
\es

The following alternative interpretation of $\H_0$, due to Corlette
\cite{cor}, Donaldson \cite{don}, Hitchin \cite{h}, and Simpson
\cite{simp}, will only be used in \S\ref{divisible}.

\bs{Theorem}
\label{three}
For $n=0$, the moduli space $\H_0$ is diffeomorphic to the moduli
space $\H$ of $\GL{2,\C}$-connections of constant central curvature $i
\omega I$, where $\omega$ is a volume form on $\rs$.  That is, $\H_0
\simeq \mu^{-1}(-I)/\GL{2,\C}$, where $\mu: \GL{2,\C}^{2g} \to
\GL{2,\C}$ is given by
$\mu(A_j,B_j)
= \prod_{j=1}^g A^{\vphantom{-1}}_j B^{\vphantom{-1}}_j A_j^{-1}
B_j^{-1}$, and $\GL{2,\C}$ acts on $\mu^{-1}(-I)$ by simultaneous
conjugation.  The natural determinant maps and universal families
coincide under this diffeomorphism.  \es

Let $\Xi$ be a fixed holomorphic line bundle over $\rs$ of degree 1,
and let $\M_n \subset \H_n$ be the subspace consisting of those $(E,
\phi) \in \H_n$ such that $\La^2 E \cong \Xi$ and $\tr \phi = 0$.  In
the case $n=0$, this is the moduli space studied by Hitchin \cite{h}.
The discussion so far, and the previous paper of the authors
\cite{ht}, refer to $\H_n$, but the remainder of this paper will
actually work with $\M_n$.  This gives equivalent information for the
following reason.  The group $\Sig = \Z_2^{2g} \subset \Jac \rs$ of
line bundles with structure group $\Z_2$ acts on $\M_n$ by tensor
product, and indeed $\H_n = (\M_n \times T^* \Jac \rs) \, / \, \Sig$.
As seen in \S1 \tag of our previous paper \cite{ht}, $H^*(\H_n) =
H^*(\M_n)^\Sig \otimes H^*(\Jac \rs)$ as rings.  To describe
$H^*(\H_n)$, therefore, it suffices to describe $H^*(\M_n)^\Sig$.
This will be the purpose of the paper.  The part of $H^*(\M_n)$ not
invariant under $\Sig$ is ignored here, but it is completely described
in a forthcoming work \cite{thad3}.

The main result of our paper on the generators \cite{ht} is the following.

\bs{Theorem}
\label{four}
The rational cohomology ring $H^*(\M_n)^\Sig$ is
generated by the universal classes, that is, the K\"unneth components
$\al_2$, $\be_2$, and $\psi_{2,j}$ of $\bar{c}_2(\Pj \bee) =
\ratio{1}{4}\,c_2 (\End \bee)$.  \es

This result has been extended to higher rank Higgs bundles by Markman
\cite{mark}.  

Following the conventions established by Newstead \cite{new}, we will
let $\al = \half \,\,\al_2$, $\be = -\ratio{1}{4} \,\,\be_2$, and
$\psi_j = \psi_{2,j}$ for $j = 1, \dots, 2g$, so that
\beq
\label{bb}
c_2 (\End \bee) = 2 \al \si - \be + 4 \sum_{j=1}^g \psi_j e_j,
\eeq
where $e_1, \dots e_{2g}$ is the usual basis for $H^1(\rs)$, and $\si
\in H^2(\rs)$ is the positive generator.

\bit{Statement of the main result}
\label{state}

Our task, then, is to give a complete set of relations between the
generators $\al$, $\be$, and $\psi_1, \dots, \psi_{2g}$.  To do so, we
must first say a little about the action of the symplectic group on
$H^*(\M_n)^\Sig$.

The group of
orientation-preserving diffeomorphisms of $\rs$
acts on $H^*(\rs)$ by automorphisms, so it
has the automorphism group of $H^*(\rs)$, namely $\Ga =
\Sp{2g, \Z}$, as a quotient \cite[p.\ 178]{mks}.

\bs{Lemma}
There is a natural action of $\Ga$ on $H^*(\M_n)^\Sig$ which fixes $\al$
and $\be$ but acts on the $\psi_j$ as the standard representation.
\es

\pf. In the case of $\M_0$, this follows immediately from \re{three},
but to extend it to $\M_n$ we describe another argument.

Let $f: \rs \to \rs$ be any orientation-preserving diffeomorphism.
The complex structure typically is not preserved by $f$, so pulling it
back induces a new complex structure $\rs'$ on the same underlying
surface.  Since Teichm\"uller space, or the moduli space of curves, is
connected, there is a path connecting $\rs$ to $\rs'$.  The
construction of the moduli space $\M_n$, and of the universal pair
$(\bee, \bphi)$, can be carried out simultaneously over all the
Riemann surfaces in this path.  Hence by homotopy invariance there is
a topological isomorphism $\bee \cong \bee'$, where $(\bee',\bphi')$
is a universal pair on $\rs'$.  The homotopy class of the isomorphism
depends only on the isotopy class of $f$.

On the other hand, if $\hat{f}: \M_n' \to \M_n$ is the map of moduli
spaces induced by $f: \rs' \to \rs$, then $(\hat{f} \times f)^*(\bee,
\bphi)$ is a universal pair over $\rs'$, and so by the uniqueness in
\re{one},
$$(\hat{f} \times f)^* \End \bee \cong \End \bee' \cong \End \bee.$$
Hence $(\hat{f} \times f)^* c_2(\End \bee) = c_2(\End \bee)$, so
$\hat{f}^* \al = \al$, $\hat{f}^* \be = \be$, and $(\hat{f}\times f)^*
\sum_j \psi_j e_j = \sum_j \psi_j e_j$.  The action of the
diffeomorphism group on $H^3(\M_n) = \langle \psi_j \rangle$ is
therefore dual to its action on $H^1(\rs) = \langle e_j \rangle$; this
factors through the standard representation of $\Ga$, which is
self-dual.  Moreover, by \re{four}, the action of the diffeomorphism
group on all of $H^*(\M_n)^\Sig$ factors through $\Ga$. \fp

The exterior square of the standard representation of $\ga$ has an
invariant element, the symplectic form.  So $\ga = -2 \sum_{j=1}^g
\psi_j \psi_{j+g} \in H^6(\M_n)$ is a $\Ga$-invariant element.  Since
the powers of the symplectic form are the only invariant elements of
exterior powers of the standard representation, we deduce the
following from \re{four}.

\bs{Proposition}
The $\Ga$-invariant part of $H^*(\M_n)^\Sig$ is generated by $\al \in
H^2$, $\be \in H^4$, and $\ga \in H^6$.
\es

Like the exterior square discussed above, the higher exterior powers
of the standard representation of $\Ga$ are reducible.  Indeed, let
$\La^k (\psi)$ be the $k$th exterior power of the standard
representation, with basis $\psi_1, \dots, \psi_{2g}$.  Define the
{\em primitive part} $\La^k_0(\psi)$ to be the kernel of the natural
map $\La^k \to \La^{2g+2-k}$ given by the wedge product with
$\ga^{g+1-k}$.  The primitive part is complementary to $\ga \La^{k-2}
\subset \La^k$, and is an irreducible representation of $\Ga$: this is
well-known for $\Sp{2g, \C}$, and so remains true for the Zariski
dense subgroup $\Ga$.  Being irreducible, it is generated by $\psi_1
\cdots \psi_k$.

For any $g, n \geq 0$, let $I^g_n$ be the ideal within the polynomial
ring $\Q[\al,\be,\ga]$ generated by $\ga^{g+1}$ and the polynomials
\beq
\label{mmm}
\rho^c_{r,s,t}
= \sum_{i=0}^{\min(c,r,s)} \, (c-i)! \,
\frac{\al^{r-i}}{(r-i)!} \,
\frac{\be^{s-i}}{(s-i)!} \,\frac{(2 \ga)^{t+i}}{i!},
\eeq
where $c = r+3s+2t-2g+2-n$,
for all $r,s,t \geq 0$ such that
\beq
\label{ppp}
r+3s+3t > 3g-3+n \mbox{\phantom{xxx} and \phantom{xxx}}
r+2s+2t \geq 2g-2+n.
\eeq
The following is then the main result of the present paper.

\bs{Theorem}
\label{ee}
As a $\Ga$-algebra,
$$H^*(\M_n)^\Sig = \bigoplus_{k=0}^g
\La^k_0(\psi)\otimes\Q[\al, \be, \ga]/I^{g-k}_{n+k}.$$
\es

The theorem enunciated in {\bf (1.2)} \tag of our previous paper
\cite{ht} is the above in the case $n=0$.  In that case the relation
of lowest degree is $\rho^g_{1,g-1,0} = g \al \be^{g-1} +
(g-1)\be^{g-2}(2\ga)$.  When $n=1$, there are two relations of lowest
degree, one of which is $\rho^{g+1}_{0,g,0}=(g+1)\be^g$.  When $n \geq
2$, the lowest degree in which a relation appears is $2(2g-2+n)$.  At
least for $s \geq r$, the relations in this degree have the
particularly simple form $\be^{s-r}(\al\be + 2 \ga)^r$.

\bit{Equivariant cohomology}
\label{equivariant}

Our main tool for studying the ring structure of $H^*(\M_n)$ is {\em
equivariant cohomology}, which we briefly review.  For a more
leisurely exposition see Atiyah-Bott \cite{ab,ab2}.

If a group --- say a Lie group --- acts on a topological space $M$,
the homotopy quotient $M_G$ is defined as the associated bundle
over the classifying space $BG$ with fiber $M$:
$$M_G = \frac{M \times EG}{G}.$$
The {\em equivariant cohomology} of $M$ is defined to be simply the
ordinary cohomology of $M_G$:
$$H^*_G(X) = H^*(X_G).$$
This is a module over $H^*(BG)$.  Restricting to any fiber gives a
natural map $H^*_G(M) \to H^*(M)$.  Note also that if $G$ acts
trivially on $M$, then $H^*_G(M) = H^*(M) \otimes H^*(BG)$.

If the action of $G$ lifts to a linear action on a vector bundle $E$,
then a vector bundle $E_G$ over $M_G$ can be defined in the obvious
way.  Thus a vector bundle equipped with such a lifting possesses
well-defined {\em equivariant characteristic classes} lying in
$H^*_G(M)$.

In our case the group acting is $\T = \cx$, so that $B\T = \C
\Pj^\infty$, and $H^*(B\T) = \Q[u]$ where $u$ is a class of degree 2.
Kirwan \cite{kir,kir2} proved the following fundamental results on
$\cx$-actions.

\bs{Theorem (Kirwan)}
\label{x}
When $\T=\cx$ acts algebraically on
a smooth quasi-projective $M$ so that $\lim_{\la \to 0} \la \cdot x$
exists for every $x$, then
\newcounter{ctr}
\begin{list}{{\rm (\roman{ctr})}}{\usecounter{ctr}}
\item there is an additive isomorphism
$$ H^i(M) \cong \bigoplus_d H^{i+r_d}(F_d),$$
where $\F = \bigcup_d F_d$ is the decomposition of the fixed-point set
into components and $r_d$ is the dimension of the subbundle of
$TM|_{F_d}$ acted on with negative weight by $\T$;
\item the Leray sequence of $M_\T \to B\T$ degenerates, so that
  $H^*_\T(M) \cong H^*(M) \otimes H^*(B\T)$ additively, and the ring
  homomorphism $H^*_\T(M) \to H^*(M)$ is surjective;
\item the restriction to the fixed-point set
$$H^*_\T(M) \lrow H^*_\T(\F) = H^*(\F)[u]$$
is an injective ring homomorphism.
\end{list}
\es

Statement (i) is perhaps most familiar in a symplectic context as
stating that the moment map is a perfect Bott-Morse function.  But
statement (iii) is equally crucial for us since it respects the ring
structure.  Together with (ii), it will tell us that a polynomial in
$\al$, $\be$, and $\psi_j$ is a relation on $\M_n$ if and only if it
is the value at $u=0$ of a polynomial in $\al$, $\be$, $\psi_j$, and
$u$ --- the {\em equivariant extension} of the relation --- whose
restriction to $H^*_G(F_d) = H^*(F_d)[u]$ is a relation for each $d$.

\bit{Fixed points of the circle action}
\label{fixed}

We will study the action of $\T = \cx$ on $\M_n$ given simply by $\la
\cdot (E, \phi) = (E, \la \phi)$.  By \re{two} this lifts to the
universal bundle, and hence the universal classes extend to
equivariant classes, which by abuse of notation, we continue to denote
$\al$, $\be$, and $\psi_j$.  They are canonical by the uniqueness in
\re{one}.

In light of \re{x}, it is vital to determine the fixed-point set for
this action.  As discussed in {\bf (10.5)} \tag of our previous paper
\cite{ht}, this would be somewhat tricky in arbitrary rank.  But now
that we have restricted attention to rank 2 (and fixed determinant),
it is not so hard.  The lemma below is proved by Hitchin \cite[7.1]{h}
for $\M_0$, but his proof generalizes directly to $\M_n$.

\bs{Lemma (Hitchin)}
\label{aa}
The components of the fixed-point set $\F$ for the $\T$-action on $\M_n$
are as follows.
\begin{list}{\rm (\roman{ctr})}{\usecounter{ctr}}
\item A component $F_0$ isomorphic to $\N$, the moduli space of stable
  bundles $E$ with $\La^2 E \cong \Xi$.  It parametrizes Higgs bundles
  of the form $(E,0)$.
\item Components $F_1, \dots, F_{g+\left[\frac{n-1}{2}\right]}$ which are fibered products
$$F_d = \rs_{2g+n-1-2d} \times_{\Jac^{2d}\rs} \Jac^d \rs,$$
where the maps $\rs_{2g+n-1-2d} \to \Jac^{2d}\rs$ and $\Jac^d \rs \to
\Jac^{2d}\rs$ are given by $D \mapsto K \Xi (n) (-D)$ and $L \mapsto
L^2$ respectively.  These parametrize Higgs bundles $(E,\phi)$ of the
form $E = L \oplus \Xi L^{-1}$,
$\phi =
\left( {\def\arraystretch{0.8} \begin{array}{cc}0&0\\s&0 \end{array}}
\right),$
where $s$ is the section of $KL^{-2} \Xi(n)$ vanishing at $D$.
\end{list}
\es

Hitchin went on to compute the cohomology of the fixed components of
type (ii) as follows.  By the Leray sequence
\beq
\label{e}
H^*(F_d) = \bigoplus_{i \in \Sig} H^*(\rs_{2g+n-1-2d}, \cl_i), \eeq
where the right-hand side consists of cohomology with local
coefficients, and $\cl_i$ runs over the flat line bundles with
structure group $\Z_2$ pulled back from $\Jac^{2d}\rs$.  If $\cl_i$ is
the trivial bundle this is simply the ordinary cohomology
$H^*(\rs_{2g+n-1-2d})$.  Otherwise, $H^k(\rs_{2g+n-1-2d}, \cl_i) =
\La^k H^1(\rs,L_i)$ if $k = 2g+n-1-2d$, and 0 if not.  Here $L_i$ runs
over the flat line bundles with structure group $\Z_2$ pulled back
from $\Jac^1 \rs$ to $\rs$ by the Abel-Jacobi map.  Hitchin shows that
for $L_i$ non-trivial, $H^0(\rs, L_i) = H^2(\rs, L_i) = 0$, and
$H^1(\rs, L_i)$ has dimension $2g-2$.

The action of $\Sig$ on $\M_n$ commutes with the $\T$-action, and
induces the trivial action on $H^*(\N)$ since it is generated by
universal classes \cite{new}.  It acts on the remaining $F_d$ as the
Galois group of the unbranched cover $F_d \to \rs_{2g+n-1-2d}$, and
the splitting \re{e} is exactly the decomposition of the cohomology
into weight spaces.

Consequently, %in the isomorphism
%$$H^*(\M_n) \cong \bigoplus_{j=0}^{g-1+n} H^*(F_d)$$
%induced by \re{x}(i), the $\Sig$-invariant part is
the $\Sig$-invariant part of $H^*(\F)$ is
$$H^*(\N) \oplus
\bigoplus_{d=1}^{g+\left[\frac{n-1}{2}\right]} H^*(\rs_{2g+n-1-2d}).$$

\bit{Symmetric products of a curve}
\label{symmetric}

The symmetric products of the curve $\rs$ thus enter into our
considerations.  So let us review some facts about the cohomology of
such a symmetric product.  Good references are the paper of Macdonald
\cite{mac} or the book of Arbarello et al.\ \cite{acgh}.

In $\rs_m \times \rs$, there is a {\em universal divisor} $\Delta$
such that $\Delta \cap (\{ D \} \times \rs) = D$.  Write its
Poincar\'e dual in terms of K\"unneth components as
$$m \si + \eta + \sum_{j=1}^{2g} \xi_j e_j,$$
where $\si$ and $e_1, \dots, e_{2g}$ are generators of $H^2(\rs)$ and
$H^1(\rs)$ respectively, so that $\eta \in H^2(\rs_m)$ and $\xi_1,
\dots, \xi_{2g} \in H^1(\rs_m)$.  A theorem of Macdonald \cite{mac}
asserts that the cohomology ring $H^*(\rs_m)$ is generated by $\eta$
and the $\xi_j$.  It is convenient to introduce $\theta_j = \xi_j
\xi_{j+g}$ and $\theta = \sum_{j=1}^g \theta_j \in H^2(\rs)$.

The group of orientation-preserving diffeomorphisms of $\rs$ acts on
$\rs_m \times \rs$.  It preserves $\Delta$, and hence the K\"unneth
components of its Poincar\'e dual.  Hence it leaves $\eta$ invariant.
Moreover, its action on the linear span of the $\xi_j$, which is
$H^1(\rs_m)$, is dual to its action on $H^1(\rs)$, and hence
factors through the quotient $\Ga = \Sp{2g,\Z}$.  There is therefore a
surjective homomorphism of $\Ga$-algebras $\La^*(\xi)[\eta] \to
H^*(\rs_m)$.  Here $\La^*(\xi)$ denotes the exterior algebra of the
standard $2g$-dimensional representation of $\Ga$, with basis vectors
$\xi_1, \dots, \xi_{2g}$.

The class $\theta$ represents the symplectic form.  So in terms of the
primitive parts $\La^k_0$ introduced in \S\ref{state}, the surjective
homomorphism above is better expressed as
$$\bigoplus_{k=0}^g \La^k_0 (\xi) \otimes \Q[\eta, \theta]
\lrow H^*(\rs_m).$$
In particular, the $\Ga$-invariant
part of $H^*(\rs_m)$ is generated by $\eta$ and $\theta$.

The following result on $H^*(\rs_m)$ will be of key importance for us.
Note that we use the term {\em total degree} to mean half the
ordinary degree of a cohomology class.

\bs{Lemma}
\label{nn}
Let $\elll$, $m$, $p$, and $q$ be non-negative integers.  If $m-g+q
\leq \elll$ and $g+p-q < \elll$, then
$$\Left( \frac{\eta^p \exp\theta}{(1+\eta)^q} \Right)_\elll = 0$$
in $H^*(\rs_m)$, where the subscript $\elll$ denotes the part in
total degree $\elll$.
\es

\pf.  Since the cup product is a homomorphism of $\Ga$-modules,
Poincar\'e duality holds for the $\Ga$-invariant part.  It therefore
suffices to check that the product of this expression with any
monomial in $\eta$ and $\theta$ evaluates to 0 on the fundamental
class of $\rs_m$.

It follows from Macdonald's results \cite{mac} that
any monomial $\eta^v \prod_j \xi_j^{w_j}$
of total degree $m$ evaluates on the fundamental class of $\rs_m$ to
1 if $w_j = w_{j+g} \leq 1$ for each $j \leq g$, and 0 otherwise.
As pointed out by Zagier \cite{pair}, this implies that for any formal
power series $A(x)$ and $B(x)$,
$$A(\eta) \exp(\theta B(\eta)) [\rs_m] = \Res_{\eta = 0}
\frac{A(\eta)(1 + \eta B(\eta))^g \, d \eta}{\eta^{m+1}}.$$
We
multiply our expression by the generating function $\exp(s
\theta)/(1+t \eta)$ for the monomials in $\eta$ and $\theta$, and ask
the coefficient of $s^i t^j$ to vanish whenever $i+j = m- \elll$:
\beqas \Coeff_{s^i t^j} \frac{\eta^p
  \exp((s+1)\theta)}{(1+\eta)^q(1+t\eta)} [\rs_m] & = & \Coeff_{s^i
  t^j} \Res_{\eta = 0} \frac{\eta^p (1 + \eta + s \eta)^g \, d
  \eta}{(1+ \eta)^q (1+t \eta)
  \eta^{m+1}} \\
& = & \const \Res_{\eta = 0} \frac{\eta^p (1+\eta)^{g-i} \eta^i \eta^j
  \, d \eta}{(1+
  \eta)^q\eta^{m+1}} \\
& = & \const \Res_{\eta = 0} \eta^{i+j+p - m-1}(1+\eta)^{g-i-q} \, d
\eta.  \eeqas
Now since $g-i-q \geq g-q-(m-\elll) \geq 0$ by hypothesis, the second
factor is a polynomial of degree $g-i-q$.  All terms therefore have
degree at most
$$(i+j+p-m-1) + (g-i-q) \leq p + (m - \elll) - m - 1 + g - q = p-\elll -
1 + g - q,$$
which is less than $-1$ by hypothesis.  \fp

\bit{Restriction of the universal classes to the fixed-point set}
\label{restriction}

In order to apply \re{x}(iii), we need to know how the equivariant
universal classes restrict to each component of the fixed-point set.
The lowest component $F_0 = \N$ is easy.  The universal pair over $\M
\times \rs$ restricts to a universal bundle over $\N \times \rs$, and
the $\T$-action restricts to a trivial action.  So $\al$, $\be$ and
$\psi_j$ restrict to classes on $\N$ defined in a like manner, and
bearing the same names.  The relations between these classes on $\N$
have been studied by many authors, notably Zagier \cite{zag}; we will
have occasion to use his results later.

The components $F_d$ for $d>0$ are handled by the following lemma.

\bs{Lemma}
\label{ll}
For $d > 0$, the restrictions of the universal classes to $F_d$ are
pulled back from the symmetric product $\rs_{2g+n-1-2d}$; indeed,
\begin{list}{{\rm(\alph{ctr})}}{\usecounter{ctr}}
\item $\al|_{F_d} = (2d-1)(\eta-u) + \theta$;
\item $\be|_{F_d} = (\eta-u)^2$;
\item $\psi_j|_{F_d} = \half(\eta-u)\,\xi_j$;
\item $\ga|_{F_d} = - \half (\eta -u)^2 \theta$.
\end{list}
\es

\pf.  We first construct an equivariant universal family $(\bee,
\bphi)$ of Higgs bundles over $F_d \times \rs$.  Since $\End \bee$ is
unique as an equivariant bundle by \re{one}, the universal
classes must restrict to the K\"unneth components of its second Chern
class.

Take the following three ingredients.  First, the line bundle
$K\Xi(n)$ over $\rs$.  Second, the universal divisor $\Delta \subset
\rs_{2g+n-1-2d} \times \rs$, or rather its associated line bundle
$\co(\Delta)$.  Third, any Poincar\'e line bundle $\cl$ over $\Jac^d
\rs \times \rs$.  Now pull all three back to $F_d \times \rs$.  By the
definition of the fibered product, $\cl^2$ and $K\Xi(n) (-\Delta)$ are
isomorphic when restricted to any fiber of the projection $F_d \times
\rs \to F_d$.  So by the push-pull formula, $\cl^2K^{-1}\Xi^{-1}(-n)
(\Delta)$ is the pull-back from $F_d$ of a line bundle, say $M$.  There
is then an element $s \in H^0(F_d \times \rs,
M\cl^{-2}K\Xi(n))$ vanishing precisely on the inverse image of
$\Delta$.

Let $\bee = \cl \oplus M \Xi \cl^{-1}$, and let $\bphi \in
H^0(\End \bee \otimes K(n))$ be given by
$\bphi =
\left( {\def\arraystretch{0.8} \begin{array}{cc}0&0\\s&0 \end{array}}
\right)$
%\left( \begin{array}{cc}0&0\\s&0 \end{array} \right)$$
with respect to the splitting.  Then $(\bee, \bphi)$ parametrizes the
pairs of the form described in \re{aa}(ii).  It is hence a universal
family.  Moreover, if $\T$ acts on the two factors with weights $1$
and $0$ respectively, then it acts on $\bphi$ by scalar
multiplication.  By \re{one}, then, $\End \bee$ is equivariantly
isomorphic to the restriction of its counterpart from $\M_n \times \rs$.

Since $\bee$ splits as a direct sum, $c_2(\End \bee) = - (c_1(M \Xi
\cl^{-1}) -  c_1(\cl))^2 = - c_1(M \Xi
\cl^{-2})^2 = - c_1(K^{-1}(-n)(\Delta))^2$.
For this to be correct equivariantly, we must include
the weights of the $\T$-action,
so the equivariant $c_1$ is the non-equivariant $c_1$
minus $u$.  It is well-known that
$$c_1(\co(\Delta)) = (2g+n-1-2d)\si +
\eta + \sum_{j=1}^{2g} \xi_j e_j:$$
see for example Arbarello et al.\ \cite{acgh}.
Hence $$c_2(\End \bee) = -((2-2g-n)\si + (2g+n-1-2d)\si +
\eta + \sum_{j=1}^{2g} \xi_j e_j - u)^2.$$
Using the identity $\left(\sum \xi_j
  e_j \right)^2 = -2 \theta \si$ and comparing coefficients with those
of \re{bb} yields the result.  \fp

All our weapons are now prepared, and we are ready to attack the main
result.  It is not a frontal assault, however.  Rather, we begin by
computing some relations quite different from the $\rho$-classes.

\bit{Some recursively defined polynomials in \boldmath $\al$, $\be$,
  $\ga$ and $u$}
\label{polynomials}

The $\Ga$-invariant subring is at the heart of the larger ring
containing it; its structure is the key to that of the whole.  Our
strategy will therefore be to look first for relations between $\al$,
$\be$, and $\ga$.  The method is curiously roundabout.  First, certain
complicated polynomials, defined recursively here in
\S\ref{polynomials}, are shown in \S\ref{relations} to be relations,
by writing down their equivariant extensions explicitly.  Then they
are shown to be expressible in terms of the much simpler polynomials
$\rho^c_{r,s,t}$ by a purely algebraic argument, given in
\S\ref{algebra}.  Not until \S\ref{wrap} does a dimension count
show that the $\rho$-classes must all be relations.

One relation which holds automatically in all $\M_n$ is $\ga^{g+1} =
0$.  This is simply because $\ga = -2 \sum_{j=1}^g \psi_j \psi_{j+g}$
and each $\psi_j^2 = 0$ by skew-commutativity.  It is therefore
convenient to view our polynomials as belonging to the ring $R =
\Q[\al, \be, \ga]/(\ga^{g+1})$.

Define polynomials $\xi^k_r$ by $\xi^k_r = 0$ for
$r <0$, $\xi^k_0 = 1$, and
\beq
\label{ff}
(r+1) \xi^k_{r+1}
= \al \xi^k_r + (r-2k) \be \xi^k_{r-1} + 2 \ga \xi^k_{r-2}
\eeq
for $r > 0$.
When $k=0$, these are the polynomials $\xi_r$ defined by Zagier
\cite{zag}, and his generating function for the $\xi_r$
extends readily.

\bs{Proposition}
\label{gg}
Define $F_0^k(x)=\sum_{r=0}^\infty \xi^k_r x^r \in
R[[x]]$.
Then
$$F_0^k(x)= (1-\be x^2)^{(2k-1)/2} \, e^{-2\ga x /\be} \left(
\frac{1+x\sqrt{\be}}
{1-x\sqrt{\be}}\right)^{(\al\be+2\ga)/2\be \sqrt{\be}}.$$
\es

\pf. From Proposition 4 of Zagier \cite{zag} we know that
$F_0^0$ satisfies the differential equation
$$(1-\be x^2) (F_0^0)'(x)=(\al +\be x + 2 \ga x^2)
F^0_0(x).$$
Now \re{ff} is equivalent to
$$(r+1)\xi_{r+1}^k-(r-1)\beta \xi^k_{r-1}=\al \xi^k_r + (1-2k) \be
\xi_{r-1}+2 \ga \xi_{r-2},$$
which shows that $F_0^k$ satisfies the differential equation
$$(1-\be x^2) (F_0^k)'(x)=(\al +(1-2k)\be x + 2 \ga
x^2)F^k_0(x),$$
with initial condition $F_0^k(0)=1$.
But $(1-\be x^2)^{k} F_0^0$ satisfies the same differential equation
\beqas
(1-\be x^2)\left((1-\be x^2)^{k} F_0^0 \right)^\prime & =&
(1-\be x^2)^{k}(1-\be x^2) (F_0^0)'-k\be x (1-\be
x^2)^{k}F^0_0 \\
& = & (\al +\be x + 2 \ga x^2)  (1-\be x^2)^{k} F_0^0
- 2k\be x  (1-\be x^2)^{k} F_0^0 \\
& = &(\al +(1-2k)\be x + 2 \ga x^2) (1-\be x^2)^{k} F_0^0,
\eeqas
and certainly
$(1-\be 0^2)^{k}F_0^0(0)=1$, so we conclude that $$F_0^k(x)=
(1-\be x^2)^{k} F_0^0(x).$$
Now substitute Zagier's generating function for $F_0^0$.
\fp

\bs{Proposition}
\label{jj}
The polynomial $\xi_r^k$ is a relation on $\N$ whenever $r \geq g +
2k$.
\es

\pf.  An equivalent form of \re{gg} is $\xi^k_r=\sum_{i=0}^k
\xi^0_{r-2i} (-\be)^{i}$, and Zagier shows that $\xi_r^0$ is a
relation on $\N$ for $r \geq g$. \fp

Now define an expression with one more index:
\beq
\label{pp}
\xi^k_{r,s} = \sum_{i=0}^s \bino{r-2k+s-i}{r-2k} \be^{s-i} \,
\frac{(2\ga)^i}{i!} \, \xi^k_{r-i}.  \eeq
Note that this is $0$ if $r<2k$.  Moreover, the $i$th term in the sum
vanishes when $i > r$, and also (as an element of $R$) when $i>g$.
Hence in particular $\xi^k_{2k,g+\elll} = \be^\elll \xi^k_{2k,g}$ in $R$
for all $\elll \geq 0$.

\bs{Proposition}
\label{mm}
Let $F^k(x,y)=\sum_{r,s=0}^\infty\xi^k_{r,s}x^r
y^s \in R[[x,y]]$.
Then
$$F^k(x,y)=\left((1-\be
y)^2-\be x^2\right)^{(2k-1)/2} \,e^{-2\ga x/\be} \,
\left( \frac{1+x\sqrt{\be}-\be
y}{1-x\sqrt{\be}-\be y}\right)^{(\al\be + 2 \ga)/2\be\sqrt \be}.$$
\es

\pf.  For fixed $r$ we have
\beqas
\sum_{s \geq 0}\xi^k_{r,s} y^s
& = &
\sum_{i=0}^r \frac{(2\ga y)^i}{i!} \, \xi^k_{r-i} \
\sum_{s=i}^\infty \bino{r-2k+s-i}{r-2k} (\be y)^{s-i} \\
& = &
\frac{1}{(1-\be y)^{r-2k+1}}
\sum_{i=0}^r \frac{(2\ga y)^i}{i!} \, \xi^k_{r-i}.
\eeqas
Multiplying by $x^r$ and summing over $r \geq 0$ we obtain
$$F^k(x,y) = (1-\be y)^{2k-1} e^{2 \ga x y/(1-\be y)}
F^k_0 \Left(\frac{x}{1-\be y}\Right),$$
and the desired result follows by substituting the formula given in
\re{gg}.  \fp

\bs{Theorem}
\label{hh}
For $r,s\geq 0$ we have
$$
\xi^k_{r,s}=\sum_{l=0}^s
(-1)^{s-l}\left[\bino{r+l}{r} + \bino{r+l-1}{r}\right]
\xi^k_{s-l}\xi^k_{r+s+l},
$$
where the binomial coefficient $\bino{-1}{0}$ is to be taken as 0.
\es

\pf.  Similar to  Zagier's Proof 1 of his Theorem 4. \fp

%\pf.  We follow Zagier's Proof 1 of his Theorem 4.  Denote the
%right-hand side of \re{ii} by $\hat{\xi}^k_{r,s}$.  We will prove the
%equality $\xi^k_{r,s} = \hat{\xi}^k_{r,s}$ by induction on $s$, the
%case $s=0$ being trivial.  A trivial binomial coefficient identity
%shows, independently of the definition of the classes $\xi^k_r$, that
%$\xi^k_{r,s}$ satisfies the recursion
%$$s\,\xi^k_{r,s}=(r-k+s)\,\be \, \xi^k_{r,s-1}+2\ga \, \xi_{r-1,s-1}.$$
%On the other hand, using \re{ff}, and taking appropriate care with the
%$l=0$ terms if $r=0$ or $1$, we find that
%\beqas
%\lefteqn{s\,\hat{\xi}^{k}_{r,s}+(r-k+s)\,\be\, \hat{\xi}^{k}_{r,s-1}-
%2\ga \,\hat{\xi}^{k}_{r-1,s-1}}\\
%&= &\sum_{l=0}^s
%(-1)^{s-l}\left[
%(s-l)\bino{r+l}{r}+(r+s+l)\bino{r+l-1}{r}\right]
%\xi^k_{s-l}\xi^k_{r+s+l}\\
%&&-\be\sum_{l=0}^{s-1}
%(-1)^{s-l-1}\left[
%(s-k-l)\bino{r+l}{r}+(r+s+l-k)\bino{r+l-1}{r}\right]
%\xi^k_{s-1-l}\xi^k_{r+s-1+l}\\
%&&-2 \ga \sum_{l=0}^{s-1}
%(-1)^{s-1-l}\left[
%\bino{r+l}{r}+\bino{r+l-2}{r}\right]\xi^k_{s-1-l}\xi^k_{r+s-2+l}\\
%&=&\sum_{l=0}^{s-1}(-1)^{s-l}\bino{r+l}{r} \\
%&&\{
%(s-l) \, \xi^k_{s-l}\xi^k_{r+s+l}-(r+s+l+1) \,
%\xi^k_{s-l-1}\xi^k_{r+s+l+1}
%\\
%&&
%-\be\left[ (s-k-l-1) \, \xi^k_{s-l-1}\xi^k_{r+s+l}-(r+s+l-k) \,
%\xi^k_{s-l-1}
%\xi^k_{r+s+l-1}\right]  \\
%&&
%-2\ga\left[
%\xi^k_{s-l-3}\xi^k_{r+s+l}-\xi^k_{s-l-1}\xi^k_{r+s+l-2}\right]
%\} \\
%&=& \sum^{s-1}_{l=0} (-1)^{s-l}\bino{r+l}{r}\{(\al
%\xi^k_{s-l-1}) \, \xi^k_{r+s+l}-\xi^k_{s-l-1}(\al
%\xi_{r+s+l})\}
%\\
%&=&0.
%\eeqas

\bit{Proof that the recursively defined polynomials are relations}
\label{relations}

After this algebraic preparation, we now find some relations
between $\al$, $\be$, and $\ga$ that can be expressed in terms of the
classes $\xi^k_{r,s}$ introduced above.  We make fundamental use of
Kirwan's theorem \re{x}(iii), which tells us that any polynomial
in the generators that vanishes on the equivariant cohomology of each
component of the fixed-point set must be a relation.

\bs{Theorem}
\label{uu}
For $n \geq 0$, let $p \in R[u]$ be an equivariant
relation on $\M_{n+2}$, that is, an element of the kernel of the natural
map to $H^*_\T(\M_{n+2})$.  Then $\partial p/ \partial u$ is an
equivariant relation on $\M_n$.
\es

\pf.  By Kirwan's theorem \re{x}(iii), it suffices to show that
$\partial p/ \partial u$ restricts to an equivariant relation on each
component of the fixed-point set of $\M_n$.

For $F_0 = \N$, this is obvious, since $p$ is also a relation on
$H^*_\T(\N) = H^*(\N)[u]$.

As for $F_d$ with $d > 0$, we may work in $H^*(\rs_m)$, where
$m=2g+n-1-2d$.  The relation $p$ restricts to a relation in
$H^*(\rs_{m+2})[u]$; moreover, since by \re{ll} the restrictions of
$\al$, $\be$, $\ga$ are polynomials in $\eta - u$ and $\theta$, the
relation can be expressed as a polynomial
$r(\phi,\theta,u)|_{\phi = \eta - u}$
such that $\partial r/ \partial u (\phi,\theta,u)|_{\phi = \eta - u}$
is the restriction of $\partial p/ \partial u$, which we want to
vanish.

The assignment $u \mapsto \eta - \phi$ induces an isomorphism
$H^*(\rs_{m+2})[u] \cong H^*(\rs_{m+2})[\phi]$, so $r(\phi, \theta,
\eta-\phi)$ is a relation in $H^*(\rs_{m+2})[\phi]$.  Observe now that
the derivative with respect to $\eta$ of any relation in
$H^*(\rs_{m+2})$ is a relation in $H^*(\rs_m)$.  This follows directly
from the list of relations given by Macdonald \cite[6.21]{mac}.
Therefore
$\partial / \partial \eta \lp r (\phi, \theta,
\eta-\phi)\rp$
is a relation in $H^*(\rs_m)[\phi]$, and hence
$$\left. \frac{\partial}{\partial \eta} \lp r (\phi, \theta,
\eta-\phi)\rp \right|_{\phi=\eta-u}
= \left. \frac{\partial r}{\partial u} (\phi,\theta,u)\right|_{\phi=\eta-u}$$
is a relation in $H^*(\rs_m)$, as desired.
\fp

\bs{Proposition}
\label{oo}
For $n \geq 0$, let $p \in R[u]$ be an equivariant
relation on $\M_n$.  Then $(u^2-\be)\,p$ is an equivariant relation on
$\M_{n+1}$.
\es

\pf.  Let $F_d$ be any component of the fixed-point set of $\M_{n+1}$.
We show that restricting $(u^2-\be)\,p$ to $F_d$ yields $0$. It is
clearly true on $F_0=\N$ since this is contained in $\M_n$.  For
$d>0$, the restriction of $u^2-\be$ to $F_d$ is $\eta \, (2u-\eta)$ by
\re{ll}. On the other hand $p$ restricted to $F_d \cap \M_n$ is
supposed to be zero.  But the image of the inclusion $F_d \cap \M_n
\subset F_d$ is Poincar\'e dual to $\eta$, since it is an \'etale
cover of the inclusion $\rs_{2g+n-1-2d} \subset \rs_{2g-2d+n}$.  Hence
$\eta$ times a relation on $F_d \cap \M_n$ is a relation on $F_d$.
\fp

These  results suggest that,  even if  we are  interested only  in the
relations on $\M_0$, it is useful to study $\M_n$ for all $n$.

\bs{Theorem}
\label{kk}
For $n\geq 0$ and $k=0,\dots,[n/2]$, the equivariant class
$$F^{k}(u,1)_{2g+2n} = \sum_{r=0}^{g+n} \xi^k_{r,g+n-r} u^r$$
is an equivariant relation on $\M_{n+2}$.
\es

\pf.  By Kirwan's theorem \re{x}(iii), it suffices to show that it
restricts to a relation on each component of the fixed-point set of
$\M_{n+2}$.

For the first component $F_0$, namely $\N$, this follows immediately
from \re{hh} and \re{jj}.

For the remaining components $F_d$ with $d>0$, use \re{ll} to restrict
\re{mm} to
$F_d$. This yields
$$F^{k}(u,1)_{2g+2n} = \left(
e^{\theta u} \,
\frac{(1-(\eta-u)(\eta-2u))^{d+k-1}}{(1-\eta(\eta-u))^{d-k}}
\right)_{2g+2n},$$
where the subscript, as in the past, denotes the part in total degree
$2g+2n$, that is, in ordinary degree $2(2g+2n)$.
To show that this vanishes in $H^*(\rs_{2g-2d+1+n}) \subset
H^*(F_d)$,
express it as
\beqas
\lefteqn{\left(e^{\theta u} \, \frac{(1-(\eta-u)(\eta-2u))^{d+k-1}}
{(1+\eta u)^{d-k} \left( 1-\frac{\eta^2}{1+\eta u} \right)^{d-k}
  }\right)_{2g+2n}}\\ \\
& = &
\left(
\sum_{i=1}^\infty \bino{d-k+i}{i}
\frac{\eta^{2i} e^{\theta u}}{(1 + \eta u)^{d-k+i}} \,
(1-(\eta-u)(\eta-2u))^{d+k-1}\right)_{2g+2n}.
\label{sum}
\eeqas
It follows immediately from \re{nn} that
$$\left( \frac{e^{\theta u} (\eta u)^{2i}}{(1+\eta u)^{d-k+i}}
\right)_{2(g-d+n+i+1-k)+j} = 0$$ for $j \geq 0$ (the 2 appearing in
the subscript since $\eta u$ and $\theta u$ are substituted for $\eta$
and $\theta$),
and hence that
$$\left( \frac{e^{\theta u} \eta^{2i}}{(1+\eta u)^{d-k+i}}
\right)_{2(g-d+n+1-k)+j} = 0$$
for $n \geq 0$.  Consequently each term
in the sum above vanishes in total degree $2g+2n$.  \fp

%\bs{Theorem} \label{tt}
%\vspace{-5.48ex}
\begin{list}{{\rm(\alph{ctr})}}%
{\usecounter{ctr}\setlength{\leftmargin}{4.7em}}
\item \bsitem{Theorem} \label{tt}
For even $n \geq 0$ and $k = 0, \dots, n/2$,
the equivariant class
$$\left( (2+u^2-\be)^{n/2-k}F^k(u,1)\right)_{2g+n+2k}$$
is an equivariant relation on $\M_{n+2}$ divisible by $u^{2k}$.
\item For odd $n \geq 0$ and $k = 0, \dots, (n\!-\!1)/2$,
the equivariant class
$$\left( (1+u^2-\be)(2+u^2-\be)^{(n-1)/2-k}F^k(u,1)\right)_{2g+n+2k+1}$$
is an equivariant relation on $\M_{n+2}$ divisible by $u^{2k+1}$.
\end{list}
%\es

\pf. Since $u$ is not a zero-divisor in $H^*_\T(\M_{n+2})$, to show
the expression in (a) is a relation it suffices to do the same for
the part in total degree $2g+2n$ of
\beqas
\lefteqn{u^{n-2k} \,
(2+u^2-\be)^{n/2-k} \,
F^k(u,1)} \\
& = &
\left((1+u^2-\be)^2-((1-\be)^2-\be u^2)\right)
^{n/2-k}F^k(u,1) \\
& = &
\sum_i \bino{\frac{n}{2}-k}{i}
(1+u^2-\be)^{2i} \left((1-\be)^2-\be u^2\right)^{n/2-k-i}
F^k(u,1) \\
& = &
\sum_i \bino{\frac{n}{2}-k}{i} (1+u^2-\be)^{2i}\,
F^{n/2-i}(u,1) \\
& = &
\sum_{i,j} \bino{\frac{n}{2}-k}{i} \bino{2i}{j}
(u^2-\be)^j \,
F^{n/2-i}(u,1).
\eeqas
By \re{kk}, $F^{n/2-i}(u,1)_{2g+2n-2j}$ is a relation on $\M_{n-j+2}$;
hence by \re{oo},
$$\left((u^2-\be)^j F^{n/2-i}(u,1) \right)_{2g+2n}$$
is a relation on $\M_{n+2}$.

The statement about divisibility is easy, since $F^k(u,1) = \sum
\xi^k_{r,s} u^r$ and $\xi^k_{r,s} = 0$ for $r < 2k$.

The proof of (b) is similar: first multiply by $u^{n-1-2k}$, compute as
before
\beqas
\lefteqn{u^{n-1-2k}\,
(1+u^2-\be)\,(2+u^2-\be)^{(n-1)/2-k}\,F^k(u,1)} \\
& = &
\sum_{i,j} \bino{\frac{n-1}{2}-k}{i} \bino{2i+1}{j}
(u^2-\be)^j\,
F^{(n-1)/2-i}(u,1),
\eeqas
then apply \re{oo} and \re{kk}.

As in (a), this is clearly divisible by $u^{2k}$, but the quotient is
further divisible by another factor of $u$.  This is because the
coefficient of $u^0$ in the quotient is
\beqas
\lefteqn{
\left(
(1-\be)
(2-\be)^{(n-1)/2-k}
\sum_{s=0}^\infty \xi^k_{2k,s}
\right)_{2g+n+1}
} \\
& = &
\left(
(1-\be)
(2-\be)^{(n-1)/2-k}
\sum_{s=g}^\infty \xi^k_{2k,s}
\right)_{2g+n+1} \\
& = &
\left(
(1-\be)
(2-\be)^{(n-1)/2-k} \, \xi^k_{2k,g}
\sum_{\elll=0}^\infty \be^\elll
\right)_{2g+n+1} \\
& = &
\left(
(2-\be)^{(n-1)/2-k} \xi^k_{2k,g}
\right)_{2g+n+1} \\
& = & 0
\eeqas
using $\xi^k_{2k,g+\elll} = \be^\elll \xi^k_{2k,g}$.
\fp

\newcounter{ctr2}
\bs{Theorem}
\label{vv}
For $n \geq -2$ even (resp. odd), $\xi^k_{r,s}$ (resp. $\xi^k_{r,s} -
\be \xi^k_{r,s-1}$) is a relation in the ordinary cohomology of
$\M_{n+2}$
\begin{list}{{\rm(\roman{ctr2})}}{\usecounter{ctr2}}
\item for $k=[n/2] - i$, $r=n-2i$, and $s=g+i$,
where $i=0,\dots, [n/2]$;
\item for $k=[n/2]+j$, $r=n+3j$, and $s=g-j$, where $j=1, \dots, g$.
\end{list}
\es

\pf.  Suppose first that $n$ is even.  For (i), just take the formula
from \re{tt}(a) with $k=n/2-i$, plug in $F^k(u,1) = \sum \xi^k_{r,s}
u^r$, divide by $u^{2k}$ and set $u=0$.  Then compute using the
definition of $\xi^k_{r,s}$, the binomial theorem, and $\ga^{g+1} =
0$.  For (ii), take the same formula on $\M_{n+2j}$ with $k = n/2 +
j$, apply \re{uu} $j$ times, and proceed as in (i).  Now suppose that
$n$ is odd.  It suffices to prove the same statement for the class
$\xi^k_{r,s} - \be \xi^k_{r,s-1} + \xi^k_{r-2,g+1}$, because in case
(i) the last term vanishes altogether, and in case (ii) it was shown
in the even case to be a relation on $\M_{n+3} \supset \M_{n+2}$.
Then everything is similar to the even case. \fp

\bit{Expressing the \boldmath $\xi$-classes in terms of
 the \boldmath $\rho$-classes}
\label{algebra}

We now have many relations on $\M_n$.  We cannot yet show that the
simple polynomials $\rho^c_{r,s,t}$ of the main theorem are relations,
but at least we can show that the relations we do have are linear
combinations of them.  Hence the goal of this section is to prove the
following purely algebraic result.

\bs{Theorem} \label{ww}
\mbox{\rm (a)} For $r \geq 2k$, $\xi^k_{r,s}$ is a linear combination
of those $\rho^{r-2k+v-w}_{u,v,w}$ with $w \leq r-2k$ and $u+3w \leq
r$.
\mbox{\rm (b)} For $r \geq 2k+1$, $\xi^k_{r,s}-\be\xi^k_{r,s-1}$ is a
linear combination of those $\rho^{r-2k+1+v-w}_{u,v,w}$ with $w \leq
r-2k+1$ and $u+3w \leq r$.  \es

It is an easy matter to check, using high-school algebra and the
equality of total degrees $r+2s+3t=u+2v+3w$, that when $n$, $k$, $r$
and $s$ are as in \re{vv}, the conditions \re{ppp} of membership
in $\ci^g_n$ are satisfied by the $\rho$-classes named above.  Hence
the relations of \re{vv} belong to $\ci^g_n$.

The proof of \re{ww} will use a generating function for the $\xi^k_r$
which generalizes a formula for the $\xi_r$ stated without proof in
Zagier's paper \cite{zag}.  Zagier kindly communicated a proof to us,
and it goes through almost verbatim for the generalization.

\bs{Lemma (Zagier)}
\label{ss}
If
$$\phi^k_m(r,p) =
\Coeff_{x^m} \, \frac{1}{\cosh^{2k} \sqrt{3x}} \,
\frac{\sqrt{3x}}{\sinh \sqrt{3x}}
\Left( \frac{\sqrt{3x}}{\tanh \sqrt{3x}} \Right)^r
\Left( \frac{1}{x} - \frac{\tanh \sqrt{3x}}{x \sqrt{3x}}\Right)^p,$$
then for $r \geq 0$,
$$\xi^k_r = \sum_{m,p} \frac{\phi^k_m(r,p)}{3^{m+p} (r-2m-3p)! \,p!}\;
\al^{r-2m-3p} \be^{m} (2 \ga)^p.$$
\es

\pf.  The formula for $\phi^k_m(r,p)$ may be abbreviated as
$\Coeff_{x^m} A(x) B(x)^r C(x)^p$.  This directly gives a generating
function for these numbers with $r$ and $p$ fixed, but to compute
$$F^k_0(t) = \sum_r \xi^k_r t^r$$
we need instead a generating function for $\phi^k_m(\elll+2m+3p,p)$
with $\elll$ and $p$ fixed and $m$ variable.  The passage from one to the
other, as usual, is by residue calculus: write $\phi^k_m(r,p)$ as
$\Res_{x=0} (A(x) B(x)^r C(x)^p dx/x^{m+1})$ and change
variables to $y =
x/B(x)^2 = (1/3) \tanh^2(\sqrt{3x})$ to get $$\phi^k_m(\elll+2m+3p,p) =
\Res_{y=0}(a(y) b(y)^\elll c(y)^p dy/y^{m+1})$$ with $a(y) = A(x(y))
x'(y)/B(y)^2$, $b(y) = B(x(y))$, $c(y) =
C(x(y)) B(x(y))^3$.
In other words,
$$\sum_m \phi^k_m(\elll+2m+3p,p) \, y^m = a(y) b(y)^\elll c(y)^p.$$
Then we need to verify
\begin{eqnarray*}F^k_0(t) & = &
\sum_{\elll,m,p \ge 0} \phi^k_m(\elll+2m+3p,p) \,
\frac{(\alpha t)^\elll}{\elll!}
\, (\beta t^2 /3)^m \, \frac{(2 \gamma t^3/3)^p}{p!} \\
& = & \sum_{\elll,p \ge 0} a(y) \,
\frac{(\alpha t b(y))^\elll}{\elll!} \, \frac{(\gamma t^3 c(y)/3)^p}{p!}  \\
& = & a(y) \exp\lp\alpha t b(y) + \gamma
t^3 c(y)/3\rp
\end{eqnarray*}
with $y = \beta t^2 /3$.
Substituting for $a$, $b$, and $c$ the formulas above,
we find $$F^k_0(t) = \cosh^{1-2k}(\sqrt{3x})
\, \exp\lp(\al \be + 2\ga) \sqrt{3x/\beta^3} -
2\gamma \tanh(\sqrt{3x})/\beta^{3/2}\rp$$ which, since the new variable $t$
is related to the original variable $x$ by $t = \sqrt{3 y/\beta} =
\beta^{-1/2} \tanh(\sqrt{3x})$, is equivalent to \re{gg}. \fp

\pf\/ of \re{ww}.  Consider first part (a).  Regarded as a polynomial
in $\al$ and $\ga$ only, each $\rho^{r-2k+v-w}_{u,v,w}$ is homogeneous
of degree $u+w$.  So let us decompose $\xi^k_{r,s}$ likewise into its
homogeneous summands relative to this $\al,\ga$-grading.  They are
nonzero only in $\al,\ga$-degree $r-2m$ for $m \geq 0$.  Indeed, using
\re{ss} and \re{pp}, we find that the part of $\xi^k_{r,s}$ in
$\al,\ga$-degree $r-2m$ equals
$$\frac{1}{3^m} \;
\sum_{i,\jay}\;%^{\min(r,s,g-t?)}
\frac{(r-2k+s-i)!}{(r-2k)!}  \;\,
\frac{\al^{r-2m-\jay}}{(r-2m-\jay)!}\;\,
\frac{\be^{s+m-\jay}}{(s-i)!} \;\,
\frac{(2\ga)^\jay}{i!\,(\jay-i)!} \;\phi^k_{m-\jay+i}(r-i,\jay-i),
$$
where we adopt the convention of summing over those indices where
the factorials all have non-negative arguments.

The $\rho$-classes having total degree $u+2v+3w=r+2s$ and
$\al,\ga$-degree $u+w=r-2m$ are of the form
$\rho^{r'+s'-2\well}_{r-2m-\well,s'-\well,\well}$ for $\well = 0,
\dots, \min([(r'+s')/2],s',r-2m)$, where we have introduced the
abbreviations $r'=r-2k$ and $s'=s+m$.
Using their definition \re{mmm}, we
may express any linear combination of these $\rho$-classes as
\beqas
&& \frac{1}{3^m} \;
\sum_{i,\well}\;%^{\min(r-2m,s')}
%\sum_{i=0}^{\min(r-2m-\well,s'-\well)}
a_\well\;
(q-2\well-i)! \;
\frac{\al^{r-i-2m-\well}}{(r-i-2m-\well)!}\;
\frac{\be^{s'-i-\well}}{(s'-i-\well)!}\;
\frac{(2\ga)^{i+\well}}{i!} \\
& = &
\frac{1}{3^m} \;
\sum_{\jay,\well}\;%^{\min(r',s+k)}
a_\well\; (q-\well-j)! \;\;\, \frac{\al^{r-2m-\jay}}{(r-2m-\jay)!}\;\;
\frac{\be^{s'-\jay}}{(s'-\jay)!}\;\;
\frac{(2\ga)^\jay}{(\jay-\well)!},
\eeqas
where $a_\well$ are
arbitrary scalars, $q = r'+s'$, and the factor of $1/3^m$ is
inserted for convenience.

At least when $s$ is large enough that $s'$ and $[(r'+s')/2] \geq
r-2m$, these span all the polynomials in $\al,\be,\ga$ of the given
total degree and $\al,\ga$-degree.  It is therefore possible to write
the part of $\xi^k_{r,s}$ in $\al,\ga$-degree $r-2m$ as a linear
combination of this kind.  The goal is to show that $a_\well=0$
whenever either $\well > r'$ or $u+3\well>r$, that is, $\well > m$.

The reader may worry that this will only prove the desired result for
$s$ large compared to $r$.  But, according to \re{mmm} and \re{pp},
the coefficient, in all of the polynomials we are concerned with, of
the monomial $\al^a\be^{s-b}\ga^c$ for fixed $a,b,c$ is a rational
function of $s$.  So if a linear dependence between these polynomials
can be established for sufficiently large $s$, then it holds for all
$s$.

To determine the scalars $a_\well$, set the coefficients of
$\al^{r-2m-\jay} \be^{s'-\jay} (2\ga)^\jay$ in the last two equations
to be equal:
$$\sum_i\;%{i=\max(0,\jay-m)}{\min(\jay,s,g?)}
\frac{(r'+s-i)!\;(s'-\jay)!}{r'!\;(s-i)!\;i!\;(\jay-i)!} \;
\phi^k_{m-\jay+i}(r-i,\jay-i)
\; = \;
\sum_\well
\frac{a_\well \;(q-\well-\jay)!}{(\jay-\well)!}.$$

Let $b_j$ be the left-hand side, and $L$ the lower triangular matrix
whose $(\jay,\well)$ entry is $(q-\well-\jay)!/(\jay-\well)!$.
Here $\jay,\well$ index the rows and columns, and run from $0$ to
$r-2m$.  In vector notation, the equation above then says $(b_\jay) =
L (a_\well)$.

The
inverse of $L$ is the lower triangular matrix whose $(\well,\jay)$
entry is $$(-1)^{\well+\jay}\frac{(q+1-2\well)}{(\well-\jay)! \:
(q+1-\well-\jay)!}.$$  Indeed, the sum that needs to be demonstrated is
$$\sum_{j=w}^{w'}(-1)^{w'+j}
\frac{(q+1-2w')(q-w-j)!}{(w'-j)!(q+1-w'-j)!(j-w)!} = \delta_{w,w'}.$$
This is obvious for $w \geq w'$.  For $w < w'$, if the summand is
denoted $N_j$, then as Shalosh B. Ekhad has pointed out \cite{ek},
$$(q+1-w-w')(w'-w)N_j = (j-w)(q+1-w-j)N_j - (j+1-w)(q-w-j)N_{j+1};$$
since the coefficient on the left is independent of $j$, and is
nonzero for $s$ and hence $q$ large, the sum telescopes.

Hence
$$a_\well \; = \; \sum_{i,\jay}
\frac{(-1)^{\well-\jay} (q+1-2\well)\; (r'+s-i)! \;
  (s'-\jay)!}{(\well-\jay)!\;(q+1-\well-\jay)!\;r'!\;(s-i)!\;i!\;(\jay-i)!}
\; \phi^k_{m-j+i}(r-i,\jay-i).$$

Now sum over all variables, and group the factorials together as
binomial coefficients, to create the grand generating function
\beqas
\lefteqn{\sum_{m,q,s,\well} \; a_\well \;\frac{(q+1-\well)!}{(q+1-2\well)\;(s'-\well)!} \;
 M^m \: Q^q \: S^s \: W^\well } \\
& = & \sum_{i,\jay,m,q,s,\well} (-1)^{\well-\jay}
\bino{r'+s-i}{r'} \bino{s'-\jay}{\well-\jay}
\bino{q+1-\well}{\jay} \bino{\jay}{i}  \;
\phi^k_{m-\jay+i}(r-i,\jay-i) \;
M^m \: Q^q \: S^s \: W^\well.
\eeqas
Using the binomial series,
we can successively eliminate the sums over $q$, $\well$, and $s$,
to obtain
$$\sum_{i,p,m} \bino{i+p}{i}
\frac{(1-WQ)^m \,(WQ^2)^{i+p}\,M^{p+m}\,S^i}{Q(1-Q)^{i+p+1}
(1-S(1-WQ))^{r'+1}}
\: \phi^k_m(r-i,p).$$
Substituting \re{ss} for the sum over $m$ yields
\beqas
\lefteqn{\frac{1}{Q(1-S(1-WQ))^{r'+1}}\sum_{i,p} \bino{i+p}{i}
\frac{(WQ^2)^{i+p} \, M^p \, S^i}{(1-Q)^{i+p+1}}}\\
& \cdot &
\frac{1}{\cosh^{2k}\sqrt{X}}\:\frac{\sqrt{X}}{\sinh \sqrt{X}} \:
\Left( \frac{\sqrt{X}}{\tanh \sqrt{X}} \Right)^{r-i} \:
\Left( \frac{3}{X} - \frac{3\tanh \sqrt{X}}{X \sqrt{X}} \Right)^p,
\eeqas
where $X = 3(1-WQ)M$.  Applying the binomial theorem again and
simplifying transforms this to
\begin{eqnarray} \label{ooo}
& &
\frac{1}{Q(1-S(1-WQ))^{r'+1}
\vphantom{\frac{\tanh\sqrt{X}}{\sqrt{X}}}} \:\;
\frac{1-WQ}{1-Q-WQ+WQ^2(1-S(1-WQ))\frac{\tanh\sqrt{X}}{\sqrt{X}}}\\
\nonumber & & \:\:\:\:\: \cdot \:\:
\frac{1}{\cosh^{2k}\sqrt{X}}\:\frac{\sqrt{X}}{\sinh \sqrt{X}} \:\Left(
\frac{\sqrt{X}}{\tanh \sqrt{X}} \Right)^r .
\end{eqnarray}

The goal is to show that the coefficient of $M^mQ^qS^sW^w$ vanishes in
the above for $q=r'+s'$ and $w > \min(r',m)$.  Since
$X=3(1-WQ)M$, it is equivalent to multiply the generating function
\re{ooo} by $3^m(1-WQ)^m$ and take the coefficient of $X^mQ^qS^sW^w$.
But the second line of \re{ooo} is a power series in $X$ only, so this
coefficient is a linear combination of the coefficients, for $n \leq
m$, of $X^nQ^qS^sW^w$ in
$$\frac{1}{Q(1-S(1-WQ))^{r'+1}
\vphantom{\frac{\tanh\sqrt{X}}{\sqrt{X}}}} \:\;
\frac{(1-WQ)^{m+1}}{1-Q-WQ+WQ^2(1-S(1-WQ))
\frac{\tanh\sqrt{X}}{\sqrt{X}}}.$$
We will show that these all vanish.

In fact, we may replace $\tanh \sqrt{X} /\sqrt{X}$ by $1+X$ in the
above.  For this can be undone by substituting a power series of the
form $c_1X + c_2 X^2 + \cdots$ for $X$; hence the coefficients of
$X^nQ^qS^sW^w$ in the former are linear combinations of $X^pQ^qS^sW^w$
in the latter, for $p \leq m$.

Taking coefficients of $X^p$, $S^s$, $W^w$, and $Q^{r'+s'}$ in the
resulting rational function yields
$$\sum_{i=0}^s (-1)^{w+i}
\bino{r'+s-p-i}{r'-p}
\bino{p+i}{i}
\bino{s'-p-i}{w-p-i}
\bino{r'+s'+1-w}{n+i}.$$
Let $F(s,i)$ be the $i$th term in the sum.  As pointed out by Shalosh
B. Ekhad \cite{ek}, if we define $G(s,i)$ by
$$\frac{i(r'+s+1-p-i)(s'+1-p-i)(r'+s'+2-w)(r'+s+s'+3-p-w-i)}%
{(s+1-i)(s'+1-w)(r'+s'+2-p-w-i)}
\,\,F(s,i),
$$
then by high-school algebra,
\beqas
\lefteqn{G(s,i+1)-G(s,i)= }\\&
(r'+s+1-w)(r'+s'+2-w)\,F(s,i)-
(s+1)(r'+s'+2-p-w)\,F(s+1,i),
\eeqas
and so, summing over $i$, we deduce that the sum satisfies a linear
recurrence relation in $s$:
$$(r'+s+1-w)(r'+s'+2-w)\sum_{i=0}^s F(s,i)
-(s+1)(r'+s'+2-p-w)\sum_{i=0}^{s+1}F(s+1,i)
= 0.$$

If $w > r'$ and $p \leq m$, the coefficient of the first sum in the
recurrence is $0$ for $s=w-r'-1$, and the coefficient of the second
sum is nonzero for all subsequent $s$.  Hence the sum vanishes for all
$s$ sufficiently large, namely $\geq w-r'$.

If $w > m$ and $p \leq m$, then every term in the sum is easily seen
to vanish for $s=0$, and for $s = -r'-m-1+p+w$ if this is positive.
The recurrence then implies that the sum is $0$ for all positive $s$.

This completes the proof of part (a); the proof of (b) is similar.
Because
$$\bino{r'+s-i}{r'} - \bino{r'+s-1-i}{r'} =
\bino{r'-1+s-i}{r'-1},$$
the grand generating function has $r'-1$ substituted for $r'$; hence
the same is true for all subsequent formulas.  \fp

\bit{The relations divisible by \boldmath $\ga$}
\label{divisible}

Many $\Ga$-invariant relations on $\M_n$ were computed in
\S\ref{relations}, but none of these relations were divisible by
$\ga$.  To find some $\Ga$-invariant relations that are divisible by
$\ga$, we revive the space $\M$ of flat connections, which is
diffeomorphic to $\M_0$ as described in \S2 \tag of our previous paper
\cite{ht}.  We will find a relationship between the cohomology at
genus $g$ and genus $g-1$.  Accordingly, we will write $\M^g$ to
indicate the dependence of $\M$ on the genus.  Let $G = \SL{2,\C}$,
and define $\mu_g: G^{2g} \to G$ by $\mu_g(A_j, B_j) = \prod
A^{\vphantom{-1}}_j B^{\vphantom{-1}}_j A_j^{-1} B_j^{-1}$.  Then
$\M^g = \mu_g^{-1}(-I)/G$, where $G$ acts by simultaneous conjugation.
The goal of this section is to prove the following

\bs{Theorem}
\label{qq}
Let $\rho \in \Q[\al, \be, \psi_j]$ be a relation on $\M^{g-1}$.  Then
$\psi_j \psi_{j+g}\rho$ for each $j \leq g$, and hence $\ga \rho$,
are relations on $\M^g$.
\es

The proof will involve the following lemma.

\bs{Lemma}
\label{rr}
The only critical value of $\mu_g$ is the identity matrix $I$.
\es

\pf.  The derivative $d\mu_g: \g^{2g} \to \g$ at $(A_j, B_j) \in
G^{2g}$ is easy to compute explicitly: see for example Goldman
\cite{gold} or Gunning \cite[Lemma 26]{gun}.  It is a sum of $g$
terms, the $k$th being conjugate to $(a_j,b_j) \mapsto (I - \Ad
A_k^{-1})b_k - (I - \Ad B_k^{-1})a_k$.  At a critical point, then,
all $g$ of these maps must fail to surject.

For $A \neq \pm I \in G = \SL{2,\C}$, it is easy to check by hand,
using Jordan canonical form, that the image of $(I - \Ad A^{-1}): \g
\to \g$ is a 2-dimensional subspace from which the eigenspaces of $A$
can be recovered, and is a subalgebra if and only if $A$ is not
diagonalizable.  If the maps are not surjective, then either these
2-dimensional subspaces must coincide, or one of $A_k$ or $B_k$ is
$\pm I$.  In the former case, $A_k$ and $B_k$ must have a common
eigenspace (if they are not diagonalizable) or eigenspaces (if they
are).  In any case, they must commute, and so $\mu_g(A_j, B_j) = I$.
\fp

Notice that the arguments of the last paragraph are special to
$\SL{2,\C}$; the situation for $\SL{r, \C}$ with $r > 2$ is more
complicated.

\pf\/ of \re{qq}.  Let $K = \SU{2}$, and let $L \subset G$ be the
locus of matrices of the form $U^{-1} D U$, where $U \in K$ and $D =
\diag(\la, 1/\la)$ for some positive real $\la$.  Then $L$ is a
smooth, contractible submanifold of $G$ whose tangent space at the
identity is $i$ times that of $K$.

Let $\cl$ be the intersection $\mu_g^{-1}(-I) \cap (G^{2g-2} \times
L^2)$.
This is preserved by the $K$-action, and there are inclusions
$$\mu_{g-1}^{-1}(-I) \times \{ I \} \times \{ I \} \subset
\cl \subset \mu_g^{-1}(-I).$$
Dividing by the $K$-action yields inclusions
$$\tilde{\M}^{g-1} \subset \cl/K \subset \tilde{\M}^g.$$ Here
$\tilde{\M}^g = \mu_g^{-1}(-I) / K$, which is a $G/K$-bundle over
$\M^g$.  Since $G/K$ is contractible, this is homotopy equivalent to
$\M^g$.  It is not hard to check, using the definition of the
universal classes in \S1 \tag of our previous paper \cite{ht}, that
$\al$, $\be$, and $\ga$, regarded as classes on $\tilde{\M}^g$,
restrict to their counterparts on $\tilde{\M}^{g-1}$.

It therefore suffices to prove the following two claims: that
$\tilde{\M}^{g-1} \subset \cl/K$ induces an isomorphism on cohomology;
and that $\cl/K \subset \tilde{\M}^g$ is Poincar\'e dual to $\psi_g
\psi_{2g}$.  For then $\psi_g \psi_{2g} \rho$ must be a relation on
$\tilde{\M}^g$.  The result follows by symmetry, since the action
of $\Ga$ on $\Q[\al, \be, \psi_j]$ certainly preserves the ideal of
relations, and there is an element of $\Ga$ taking $\psi_g \psi_{2g}$
to $\psi_j \psi_{j+g}$ for each $j$.

To prove the first claim, first note that $\cl$ can be regarded as the
fibered product $G^{2g-2} \times_G (L \times L)$, where the map
$G^{2g-2} \to G$ is $\mu_{g-1}$ and the map $L \times L \to G$ is
$(A,B) \to - ABA^{-1}B^{-1}$.  A direct computation shows that no two
$A,B \in L$ have $ABA^{-1}B^{-1} = -I$; hence by \re{rr} the map $L
\times L \to G$ never touches a critical value of $\mu_{g-1}$.  So
$\cl$ is locally trivial over $L \times L$ with fiber
$\mu_{g-1}^{-1}(-I)$.  Since $L \times L$ is contractible, this implies
that $\cl$ is homeomorphic to $L \times L \times \mu_{g-1}^{-1}(-I)$.

It follows that $\tilde{\M}^{g-1}$ and $\cl/K$ are homotopy equivalent.
Indeed, they are homotopy equivalent to the homotopy quotients
$(\mu_{g-1}^{-1}(-I) \times EK)/K$ and $(\cl \times EK)/K$
respectively, and the latter retracts onto the former, since $BK$ is a
direct limit of manifolds and $L \times L$ is contractible.

To prove the second claim, first note that $L$ meets $K\subset G$
transversely at the single point $I$.  It is therefore Poincar\'e dual
to the standard generator of $H^3(G,\Z)$.  We can now either imitate
the argument given by the second author \cite[Prop.\ 19.3]{thad2} for
the $\SU{2}$ space $\N^g = (\mu_g^{-1}(-I) \cap K^{2g})/K$, or simply
use that result.  It tells us that the natural maps in the top row of
the diagram
$$\begin{array}{ccccc}
K^{2g} & \longleftarrow & \mu_g^{-1}(-I) \cap K^{2g}
& \longrightarrow & \N^g \\
\down{} & & \down{} & & \down{} \\
G^{2g} & \longleftarrow & \mu_g^{-1}(-I)
& \longrightarrow & \tilde{\M}^g
\end{array}$$
induce isomorphisms on $H^3$ under which the generator of the $j$th
copy of $H^3(K, \Z)$ corresponds to $\psi_j$.  Since the outer columns
also induce isomorphisms on $H^3$, so does every map in the diagram.

Since the pull-back by inclusion is Poincar\'e dual to transverse
intersection, it now suffices to check that $G^{2g-2} \times L^2$ is
transverse to $\mu_g^{-1}(-I)$, or equivalently, that at every point
of $\cl$ the derivative $d \mu_g$ remains surjective when restricted
to the tangent space to $G^{2g-2} \times L^2$.  But this is true even
if we restrict further to the tangent space to $G^{2g-2}$, since as
stated before we are at a regular value of $\mu_{g-1}$.  \fp

\bit{The cohomology not fixed by \boldmath $\Ga$}
\label{gamma}

Everything so far has been about the $\Ga$-invariant part of
$H^*(\M_n)^\Sig$, generated by $\al$, $\be$, and $\ga$.  Now it is
time to say something about the non-invariant part and the classes
$\psi_j$.

We begin with a result relating the non-invariant parts of the
cohomology of the symmetric product $\rs_m$ at genus $g$ to the
invariant part at lower genera.

\bs{Lemma}
As a $\Ga$-module, the cohomology of the symmetric product
$\rs^g_{m\vphantom{k}}$ has the form
$$H^*(\rs^g_{m\vphantom{k}})
= \bigoplus_{k=0}^g \La^k_0(\xi) \otimes \Q[\eta,\theta]/I(\rs^{g-k}_{m-k}),$$
where $I(\rs^{g-k}_{m-k})$ is the ideal of relations between $\eta$
and $\theta$ in $\rs^{g-k}_{m-k}$, with the convention that
$I(\rs^{g-k}_{m-k}) = \Q[\eta,\theta]$ if $m-k < 0$.
\es

\pf.  As shown in \S\ref{symmetric}, there is a surjection of
$\Ga$-algebras
$$\bigoplus_{k=0}^g \La^k_0 (\xi) \otimes \Q[\eta, \theta] \lrow
H^*(\rs^g_{m\vphantom{k}}),$$
where $\La^k_0$, being irreducible, is spanned by the
orbit of $\xi_1 \xi_2 \cdots \xi_k$ under $\Ga$.  It therefore
suffices to show that a polynomial $p(\eta, \theta)$ is a relation on
$\rs^{g-k}_{m-k}$ if and only if $p(\eta, \theta) \xi_1 \xi_2 \cdots
\xi_k$ is a relation on $\rs^g_{m\vphantom{k}}$.

By Poincar\'e duality the latter is true if and only if for all
polynomials $q(\eta, \xi_j)$,
\beq
\label{f} p(\eta, \theta) \xi_1
\xi_2 \cdots \xi_k q(\eta, \xi_j) [\rs^g_{m\vphantom{k}}] = 0.  \eeq
Now it follows from the description of $H^*(\rs^g_{m\vphantom{k}})$ in
Macdonald \cite{mac} that $\prod_{j=1}^{2g} \xi_j^{p_j} \eta^q
[\rs^g_{m\vphantom{k}}] = 1$ if $\sum_{j=1}^g p_j + q = m$ (so that
the degrees match) and $p_j = p_{j+g} \leq 1$ for each $j \leq g$ (so
that it becomes a monomial in $\eta$ and $\theta_j = \xi_j
\xi_{j+g})$.  Otherwise it equals 0.  Hence in \re{f} we only need to
consider the case
$$q(\eta,\xi_j) = \xi_{g+1} \cdots \xi_{g+k} r(\eta, \theta_{k+1},
\dots, \theta_g).$$
But $q$ can be averaged with its images under all
permutations of $\theta_{k+1}, \dots, \theta_g$ without changing the
value of \re{f}.  Hence we only need to consider the case where $r$ is
a polynomial in $\eta$ and $\theta_{k+1} + \dots + \theta_g$.  But
then
\beqas
p(\eta, \theta) \xi_1 \xi_2 \cdots \xi_k q(\eta, \xi_j)
[\rs^g_{m\vphantom{k}}] & = & (-1)^k p(\eta, \theta) \theta_1 \theta_2 \cdots
\theta_k
r(\eta, \theta_{k+1} + \dots + \theta_g) [\rs^g_{m\vphantom{k}}] \\
& = & (-1)^k p(\eta, \theta_{k+1} + \dots + \theta_g) \theta_1
\theta_2 \cdots \theta_k r(\eta, \theta_{k+1} + \dots + \theta_g)
[\rs^g_{m\vphantom{k}}].
\eeqas
This always vanishes if $k > m$.  Otherwise, it
equals $(-1)^k p(\eta, \theta) r(\eta, \theta) [\rs^{g-k}_{m-k}]$.
This vanishes for all $r$ if and only if $p(\eta, \theta) q(\eta,
\xi_j) [\rs^{g-k}_{m-k}]$ vanishes for all polynomials $q$ in $\eta$
and $\xi_j$, since it does not alter the latter expression to replace
$q$ with its projection on the $\Ga$-invariant part, which is a
polynomial in $\eta$ and $\theta$.  Again by Poincar\'e duality, this
is equivalent to $p(\eta,\theta) = 0$ in $\rs^{g-k}_{m-k}$. \fp

\bs{Theorem}
\label{xx}
As a $\Ga$-module, the $\Sig$-invariant part of the $\T$-equivariant
cohomology of $\M^g_n$ has the form
$$H_\T^*(\M^g_n)^\Sig \cong \bigoplus_{k=0}^g
\La^k_0(\psi)\otimes\Q[\al, \be, \ga, u]/I_\T(\M^{g-k}_{n+k}).$$
Consequently, as a $\Ga$-module, the $\Sig$-invariant part of the
ordinary cohomology of $\M^g_n$ has the form
$$H^*(\M^g_n)^\Sig \cong \bigoplus_{k=0}^g
\La^k_0(\psi)\otimes\Q[\al, \be, \ga]/I(\M^{g-k}_{n+k}).$$
\es

\pf.  First we show that if $\la\in \La^k_0(\psi)$ and $r\in
I_\T(\M^{g-k}_{n+k})$, then $\la r \in I_\T(\M^g_n)$.  Certainly $r$
restricts to relations between $\al$, $\be$, $\ga$ and $u$ on $\N^g =
F^g_0$ and between $\eta$, $\theta$ and $u$ on the remaining fixed
components $F^g_d \subset \M^g_n$.  On the other hand, by Kirwan's
theorem \re{x}(iii) it suffices to show that $\la r$ restricts to
similar relations on the fixed components $F^{g-k}_d$ of
$\M^{g-k}_{n+k}$.

The case of $d > 0$ follows immediately from the lemma.  As for $d=0$,
Proposition 2.5 of King-Newstead \cite{kn} asserts that if $\la\in
\La^k_0(\psi)$ and $r$ is a relation between $\al$, $\be$, $\ga$ on
$\N^g$, then $\la r$ is a relation on $\N^{g-k}$.  This is exactly
what is needed.

The left-hand side is therefore a quotient of the right-hand side.  To
complete the proof, it remains only to check that the
$\Sig$-invariant, $\T$-equivariant Poincar\'e polynomials agree.  But
if $P = \sum_i t^i \dim H^i$, then
\beqas
\lefteqn{(1-t^2)P_\T^\Sig(\M^g_n)} \\
&=& P^\Sig(\M^g_n)\\
&=& P^\Sig(\N^g)
+ \sum_{d=1}^{g+\left[\frac{n-1}{2}\right]}
t^{2g+2d-2}P(\rs_{2g+n-1-2d}^g)\\
&=& \sum_{k=0}^g\left(\bino{2g}{k}-\bino{2g}{k-2}\right) \LEFT(t^{3k}
P^\Sig(\N^{g-k})+ \sum_{d=1}^{g+\left[\frac{n-1}{2}\right]}
t^{2g+2d-2} t^k P(\rs_{2g+n-1-2d-k}^{g-k})\RIGHT)\\
&=& \sum_{k=0}^g\left(\bino{2g}{k}-\bino{2g}{k-2}\right)
t^{3k}\LEFT(P^\Sig(\N^{g-k})+\sum_{d=1}^{g+\left[\frac{n-1}{2}\right]}
t^{2(g-k)+2d-2}P(\rs^{g-k}_{2(g-k)+n-1-2d+k})\RIGHT) \\
&=& \sum_{k=0}^g\left(\bino{2g}{k}-\bino{2g}{k-2}\right)
t^{3k} P^\Sig(\M^{g-k}_{n+k})\\
&=& (1-t^2)\sum_{k=0}^g\left(\bino{2g}{k}-\bino{2g}{k-2}\right) t^{3k}
P_\T^\Sig(\M^{g-k}_{n+k}),
\eeqas
using Kirwan's theorem on the Leray sequence \re{x}(ii) in steps 1 and
6, the perfection of the Morse stratification \re{x}(i) in steps 2 and
5, the lemma and the result of King-Newstead in step 3, and
high-school algebra in step 4.  \fp

\bit{Wrap-up}
\label{wrap}

At last we can show that the polynomials $\rho^c_{r,s,t}$ of \re{mmm}
are relations, using every tool at our disposal: \re{ww}, \re{qq},
\re{xx}, and a dimension count.  Recall that $\ci^g_n$ is the ideal of
$\rho$-classes introduced in \S\ref{state}.

\bs{Proposition}
Every element of $\ci^g_n$ is a relation on $\M^g_n$.
\es

\pf.  It is actually more convenient to work with $n+2$ than $n$, so
let $n \geq -2$.  Then every element of $\ci^g_{n+2}$ has total degree
$\geq 2g+n$.  In degree $2g+n$, $\ci^g_{n+2}$ is spanned by the
elements $\rho^{g+i}_{n-2i,g+i,0}$ for $i= 0,\dots,[n/2]$.  The
degrees with respect to $\al$ are all different, so these are linearly
independent.  The relations $\xi^{[n/2]-i}_{n-2i,g+i}$ of type (i)
given in \re{vv} are all of degree $2g+n$, also number $[n/2] + 1$,
and are linearly independent for the same reason.  By \re{ww} they are
in the linear span of the $\rho^{g+i}_{n-2i,g+i,0}$.  Hence this
equals the linear span of the relations of type (i), so all the
$\rho^{g+i}_{n-2i,g+i,0}$ are relations.

The diagram below is intended to help the reader visualize the ideal
of relations.  Each dot represents one of the $\rho^c_{r,s,t}$.  The
total degree is the vertical coordinate, and $r$ is the horizontal
coordinate.  The two edges of the dotted region reflect the two
constraints imposed by \re{ppp}.  To avoid having to draw a third
dimension, only those relations with $t=0$ have been shown.  The
dotted region for any fixed $t>0$ would look similar, only translated in
a northwesterly direction.  The relations established in the previous
paragraph are those in the bottom row.

The rest of the proof proceeds by induction on the total degree.  We
have already seen that the part of $\ci^g_n$ in total degree $2g+n$
consists entirely of relations.  Now fix $j>0$ and consider the part
of $\ci^g_{n+2}$ in degree $2g+n+j$.  Assume by induction that
for all $g$ and $n$, the parts of $\ci^g_{n+2}$ in degree $<2g+n+j$ are
known to be relations on $\M^g_{n+2}$.

In particular, if $\rho^c_{r,s,t} \in \ci^g_{n+2}$ has degree $2g+n+j$
and $t>0$, then $\rho^c_{r,s,0} \in \ci^{g-t}_{n+2}$ has degree $2(g-t)
+ n + (j-t)$, and hence is a relation on $\M^{g-t}_{n+2}$.  By
\re{xx}, $\psi_1 \cdots \psi_{n+2} \rho^c_{r,s,0}$ is a relation on
$\M^{g-t+n+2}_0$, so by \re{qq}, $\psi_1 \cdots \psi_{n+2}
\rho^c_{r,s,t}$ is a relation on $\M^{g+n+2}_0$, and by \re{xx} again,
$\rho^c_{r,s,t}$ is a relation on $\M^g_{n+2}$. Because these
relations have $t>0$, they are not shown in the diagram.

Only the $\rho^c_{r,s,t}$ with $t=0$ remain.  Consider first those
relations of the form $\rho^c_{0,s,0} = c!/s! \be^s \in \ci^g_{n+2}$.
These are the relations at the left-hand edge of the diagram.  If the
degree equals $2g+n+j$, then since $j>0$,
$$2s=2g+n+j \geq 2g-2+(n+3)$$
and
$$3s=3g+{\textstyle \frac{3}{2}} n +
{\textstyle \frac{3}{2}} j > 3g-3+(n+3),$$
so in fact $\rho^c_{0,s,0} \in \ci^g_{n+3}$.
Since it has degree $2g+n+j=2g+(n+1)+(j-1)$, by the induction
hypothesis again it is a relation on $\M^g_{n+3}$, and hence on
$\M^g_{n+2} \subset \M^g_{n+3}$.

\begin{figure}

\begin{center}
\begin{picture}(300,225)
\thicklines
\put(0,0){\vector(0,1){225}}
\put(-35,210){\small $r+2s$}
\put(0,0){\vector(1,0){300}}
\put(305,-2){\small $r$}
\multiput(10,60)(20,0){7}{\circle*{3}}
\multiput(0,90)(20,0){9}{\circle*{3}}
\multiput(10,120)(20,0){10}{\circle*{3}}
\multiput(0,150)(20,0){12}{\circle*{3}}
\multiput(10,180)(20,0){13}{\circle*{3}}
\multiput(0,210)(20,0){15}{\circle*{3}}

\thinlines
\multiput(12,66)(20,0){7}{\vector(1,3){6}}
\multiput(2,96)(20,0){9}{\vector(1,3){6}}
\multiput(12,126)(20,0){10}{\vector(1,3){6}}
\multiput(2,156)(20,0){12}{\vector(1,3){6}}
\multiput(12,186)(20,0){13}{\vector(1,3){6}}
\end{picture}
\end{center}

\end{figure}

What if $t=0$ but $r >0$?  It is easily checked that
$$r \rho^c_{r,s,0} = c \,\al\rho^{c-1}_{r-1,s,0} +
(c-r)\rho^{c-2}_{r-1,s-1,1}.$$
Now if $r+2s=2g+n+j$ and $1 \leq r \leq
n+3j-2$ (so that we are not at the right-hand edge of the diagram), we
know from the induction hypothesis that $\rho^{c-1}_{r-1,s,0}$ is a
relation.  And certainly $\rho^{c-2}_{r,s-1,1}$ is a relation, since
it belongs to $\ci^g_{n+2}$ and has $t>0$.  Hence under these
circumstances $\rho^c_{r,s,0}$ is a
relation.  This
fills in the interior of the diagram; the arrows depict multiplication
by $\al$ (modulo $\ga$), which takes $\rho^{c-1}_{r-1,s,0}$ to
$\rho^c_{r,s,0}$.

Only one class of total degree $2g+n+j$ remains in $\ci^g_{n+2}$.
This is $\rho^g_{n+3j,g-j}$, at the right-hand edge of the diagram.
We know from \re{vv}(ii) that the class $\xi^{[n/2]+j}_{n+3j,g-j}$ is
a relation if $n$ is odd, and that $\xi^{[n/2]+j}_{n+3j,g-j} - \be
\xi^{[n/2]+j}_{n+3j,g-j-1}$ is a relation if $n$ is even.  Since the
leading term with respect to $\al$ of $\xi^k_{r,s}$ is ${r-2k+s
\choose r-2k} \al^k/k!$, the monomial $\al^{n+3j} \be^{g-j}$ appears
in these relations with a nonzero coefficient.  On the other hand, by
\re{ww} these relations can be expressed as a linear combination of
$\rho^c_{r,s,t} \in \ci^g_{n+2}$.  Since $\rho^g_{n+3j,g-j}$ is the
class of maximal degree $n+3j$ with respect to $\al$ among these, its
coefficient in this combination must be nonzero.  It is therefore
indeed a relation.  This completes the proof.  \fp

\bs{Proposition}
Every relation between $\al$, $\be$, $\ga$ on $\M^g_n$ is an element of
$\ci^g_n$.
\es

\pf.  Since the converse has just been shown, it suffices to
show that $\dim H^I(\M^g_n) = \dim \Q[\al,\be,\ga]/\ci^g_n$, where
$H^I$ denotes the part of $H^*$ invariant under the
action of both $\Sig$ and $\Ga$.

(In the proofs of the analogous statement for $\N^g$
\cite{bar,kn,st,zag}, it has been customary to show that all the Betti
numbers match up.  This can certainly be done for $\M^g_n$, but the
cruder statement about overall dimension is clearly sufficient.)

Now since by \re{x}(i) the $\cx$ action is perfect, the dimension of
$H^I(\M^g_n)$ may be expressed as a sum over the fixed components
enumerated in \re{aa}:
$$\dim H^I(\M^g_n) = \dim H^I(\N^g) +
\sum_{d=1}^{g+\left[\frac{n-1}{2}\right]}
\dim H^I(\rs_{2g+n-1-2d}).$$
The Poincar\'e polynomial of $H^I(\N^g)$ is
$$\frac{(1-t^{2g})(1-t^{2g+2})(1-t^{2g+4})}{(1-t^2)(1-t^4)(1-t^6)};$$
this is clear, for example, from the presentation with three
generators and three relations \cite{bar,kn,st,zag}.  To find $\dim
H^I(\N^g)$ we want to substitute $t=1$.  Of course 0 appears in the
denominator, but the limit as $t \to 1$ can easily be evaluated to
${g+2 \choose 3}$ by substituting $t^2 = 1 + \epsilon$, then using
$(1+\epsilon)^k=1+k\epsilon+O(\epsilon^2)$.

As for $\dim H^I(\rs_m)$, it follows easily from the discussion in
Arbarello et al.\ \cite[VII B]{acgh} that this equals $\left[
  \frac{m+2}{2} \right] \left[ \frac{m+3}{2} \right]$ if $m \leq 2g-1$
and $(g+1)(m-g+1)$ if $m \geq 2g-1$.  Hence \beqas
\lefteqn{\dim H^I(\M^g_n)} \\
& = & \bino{g+2}{3} + \sum_{d=1}^{\left[ \frac{n-1}{2}\right]}
(g+1)(g-2d+n)
+\sum_{d=\left[\frac{n-1}{2}\right]+1}^{\left[\frac{n-1}{2}\right]+g}
\left(g-d+1+{\textstyle \left[ \frac{n-1}{2} \right]}\right)
\left(g-d+1+{\textstyle \left[ \frac{n}{2} \right]}\right).  \eeqas

On the other hand,
$$\dim \Q[\al, \be, \ga]/\ci^g_n = \sum_{r,s,t} \, 1$$
where the right-hand sum runs over all non-negative $r,s,t$ such that
$t \leq g$ and $r+3s+3t \leq 3g-3+n$ {\em or} $r+2s+2t < 2g-2+n$.
This can be re-written as
\beqas
\lefteqn{
\sum_{t=0}^{g
\vphantom{\left[\frac{n}{2}\right]}}
\sum_{s=0}^{g-1+\left[ \frac{n}{2}\right] -t}
\sum_{r=0}^{\scriptstyle \max(3g-3+n-3s-3t,
\atop \scriptstyle 2g-3+n-2s-2t)}
1
} \\
& = & \sum_{t=0}^g \left(
\sum_{s=0}^{g-1-t}(3g-2+n-3s-3t) +
\sum_{s=0}^{g-1+\left[ \frac{n}{2}\right]-t}(2g-2+n-2s-2t)
\right).
\eeqas
It is straightforward, using high-school algebra and the identities
$\sum_{d=1}^k d = k^2/2 + k/2$ and $\sum_{d=1}^k d^2 =
k^3/3+k^2/2+k/6$, to show that this equals the above formula for $\dim
H^I(\M^g_n)$.  \fp

\pf\/ of \re{ee}.
The two results above show that the ideal $I(\M^g_n)$
of $\Ga$-invariant relations
on $\M^g_n$ is precisely $\ci^g_n$.  Now apply \re{xx}. \fp

\bit{Relationship with other papers}
\label{tie-up}

The present paper is closely related to several other works by the
authors; in this final section we indicate briefly a few of the points
of contact.

The first author has constructed a compactification $\overline{\H}_n$
of the moduli space of Higgs bundles \cite{hausel1}, by adding a
divisor $Z_n$ at infinity which is a quotient by $\T = \cx$ of an open
subset of $\H_n$.  Indeed, $\overline{\H}_n$ itself is a quotient by
$\T$ of $\H_n \times \C$.  Many of the constructions given herein
apply to $\overline{\H}_n$ and $Z_n$.  In particular, there are direct
limits $\overline{\H}_\infty$ and $Z_\infty$.  Just as $\H_\infty$ is
shown in {\bf (9.7)} of our previous paper \cite{ht} to be homotopy
equivalent to $B \overline{\G}$, we expect $\overline{\H}_\infty$, and
also $Z_\infty$, to be homotopy equivalent to $B \G$, the classifying
space of the full gauge group.

The cohomology rings of $\overline{\H}_\infty$ and $Z_\infty$, and
hence those of $\overline{\H}_n$ and $Z_n$, will have generators like
those of $\H_n$, but with one additional generator $h \in H^2$,
corresponding to the class discarded in the proof of {\bf (10.1)} of
our previous paper.  Indeed, the quotient map $H^*_\T(\H_n) \to
H^*(Z_n)$ is surjective, and $h$ is the image of $u$.  It can be shown
that the kernel is the image in $H^*_\T(\H_n)$ of the compactly
supported cohomology.  The relations between the generators in
$H^*(Z_n)$ are therefore of two types: those coming from relations in
$H^*_\T(\H_n)$, and those coming from the compactly supported
cohomology.  The former are covered by the results of this paper.  As
for the latter, they ought to be determined by the results of another
paper of the first author \cite{hausel2}, in which the intersection
pairings between the cohomology and the compactly supported cohomology
of $\M_0$ were computed (and shown to vanish).

By studying the stratification of $\H_\infty$ by the HN type of the
underlying bundle $E$ --- not $(E,\phi)$ --- in the rank 2 case, the
first author has been able to extract the relations in the cohomology
of the lowest stratum, which retracts onto $\N$.  This recovers the
description of the ring $H^*(\N)$, given by several authors
\cite{bar,kn,st,zag}, in an especially simple and geometrical fashion.
It will be described in a forthcoming paper \cite{mum}, where it will
also be shown how consideration of the Mumford conjecture leads
to a natural geometrical proof of the generation theorem for the
moduli space of Higgs bundles in ranks 2 and 3.

The second author has studied moduli problems providing smooth
resolutions of the upward and downward flows from the components of
the fixed-point set of $\H_n$ in the rank 2 case.  The downward flow
is intriguing because it can be interpreted as a master space of
Bradlow pairs, but the upward flow is particularly related to the
present paper.  It is relatively simple to describe, but it contains the
parts of $H^*(\H_n)$ not invariant under $\Sig = \Z_2^{2g}$, in
the sense that none of them are killed by the restriction map.  These
upward flows can therefore be used to complete the description of the
cohomology rings $H^*(\H_n)$, by characterizing the part not invariant
under $\Sig$.  In fact this is not so difficult, since for dimension
reasons these classes have square 0, and are killed by the $\psi_j$.
So the products with $\al$ and $\be$ are all that must be computed.
Details will appear in a forthcoming paper \cite{thad3}.


\newpage

\begin{thebibliography}{99}
{\small

\bibitem{acgh}{\sc E. Arbarello, M. Cornalba, P.A. Griffiths,
{\rm and} J. Harris},
{\it Geometry of algebraic curves,} volume I,
Springer-Verlag, 1985.

\bibitem{ab}{\sc M.F. Atiyah {\rm and} R. Bott,}
The Yang-Mills equations over Riemann surfaces,
{\it Philos.\ Trans.\ Roy.\ Soc.\ London Ser.\ A} 308 (1982) 523--615.

\bibitem{ab2}{\sc M.F. Atiyah {\rm and} R. Bott,}
The moment map and equivariant cohomology,
{\it Topology} 23 (1984) 1--28.

\bibitem{bar}{\sc V. Yu.\ Baranovski\u\i},
The cohomology ring of the moduli space of stable bundles with odd
determinant,
{\it Izv.\ Ross.\ Akad.\ Nauk Ser.\ Mat.\ }58 (1994) 204--210;
translation in {\it Russian Acad.\ Sci.\ Izv.\ Math.\ }45 (1995)
207--213.

\bibitem{cor}
{\sc K. Corlette},
Flat $G$-bundles with canonical metrics,
{\it J. Differential Geom.\ }28 (1988) 361--382.

\bibitem{don}
{\sc S.K. Donaldson},
Twisted harmonic maps and the self-duality equations,
{\it Proc.\ London Math.\ Soc.\ }(3) 55 (1987) 127--131.

\bibitem{ek}
{\sc S.B. Ekhad}, A proof of a recurrence,
{\it EKHAD}, available from
{\tt http://www.math.temple. edu/\~{}zeilberg/programs.html/}.

\bibitem{gold}
{\sc W.M. Goldman},
The symplectic nature of fundamental groups of surfaces,
{\it Adv.\ Math.\ }54 (1984) 200--225.

\bibitem{gs}
{\sc D. Grayson {\rm and} M. Stillman},
{\it Macaulay 2}, available from {\tt
  http://www.math.uiuc.edu/ Macaulay2/}.

\bibitem{gun}
{\sc R.C. Gunning},
{\it Lectures on vector bundles over Riemann surfaces, }
Princeton, 1967.

\bibitem{hausel1}
{\sc T. Hausel},
Compactification of moduli of Higgs bundles,
{\it J. Reine Angew.\ Math.\ }503 (1998) 169--192.

\bibitem{hausel2}
{\sc T. Hausel},
Vanishing of intersection numbers on the moduli space of Higgs
bundles,
{\it Adv.\ Theor.\ Math.\ Phys.\ }2 (1998) 1011--1040.

\bibitem{mum}
{\sc T. Hausel},
Geometric proof of the Mumford conjecture, in preparation.

\bibitem{ht}
{\sc T. Hausel {\rm and} M. Thaddeus},
On the generators of the cohomology ring
of the moduli space of rank 2 Higgs bundles,
preprint {\tt math.AG/0003093}.

\bibitem{h}
{\sc N.J. Hitchin,}
The self-duality equations on a Riemann surface,
{\it Proc.\ London Math.\ Soc.\ }(3) 55 (1987) 59--126.

\bibitem{kn}
{\sc A.D. King {\rm and} P.E. Newstead},
On the cohomology ring of the moduli space of rank $2$ vector bundles
on a curve,
{\it Topology }37 (1998) 407--418.

\bibitem{kir}
{\sc F.C. Kirwan},
{\it Cohomology of quotients in algebraic and symplectic geometry, }
Math.\ Notes 31, Princeton, 1984.

\bibitem{kir2}
{\sc F.C. Kirwan},
The cohomology rings of moduli spaces of bundles over Riemann
surfaces,
{\it J. Amer.\ Math.\ Soc.\ }5 (1992) 853--906.

\bibitem{mks}
{\sc W. Magnus, A. Karrass, {\rm and} D. Solitar},
{\it Combinatorial group theory},
Interscience, 1966.

\bibitem{mac}
{\sc I.G. Macdonald},
Symmetric products of an algebraic curve,
{\it Topology }1 (1962) 319--343.

\bibitem{mark}
{\sc E. Markman},
Generators of the cohomology ring of moduli spaces of sheaves on
symplectic surfaces,
preprint {\tt math.AG/0009109}.

\bibitem{munoz}
{\sc V. Mu\~noz},
Ring structure of the Floer cohomology of $\Sigma\times{S}\sp 1$,
{\it Topology }38 (1999) 517--528.

\bibitem{new}
{\sc P.E. Newstead},
Characteristic classes of stable bundles
over an algebraic curve,
{\it Trans.\ Amer.\ Math.\ Soc.\ }169 (1972) 337--345.

\bibitem{nit}
{\sc N. Nitsure},
Moduli space of semistable pairs on a curve,
{\it Proc.\ London Math.\ Soc.\ }(3) 62 (1991) 275--300.

\bibitem{st}
{\sc B. Siebert {\rm and} G. Tian},
Recursive relations for the cohomology ring of moduli spaces of stable
bundles,
{\it Turkish J. Math.\  }19 (1995) 131--144.

\bibitem{simp3}
{\sc C.T. Simpson},
Higgs bundles and local systems,
{\it Inst.\ Hautes Etudes Sci.\ Publ.\ Math.\ }75 (1992) 5--95.

\bibitem{simp}
{\sc C.T. Simpson},
Moduli of representations of the fundamental group of a smooth
projective variety, I,
{\it Inst.\ Hautes Etudes Sci.\ Publ.\ Math.\ }79 (1994) 47--129.

\bibitem{simp2}
{\sc C.T. Simpson},
Moduli of representations of the fundamental group of a smooth
projective variety, II,
{\it Inst.\ Hautes Etudes Sci.\ Publ.\ Math.\ }80 (1995) 5--79.

\bibitem{cft}
{\sc M. Thaddeus},
Conformal field theory and the cohomology of the moduli space of
stable bundles,
{\it J. Differential Geom.\ }35 (1992) 131--149.

\bibitem{pair}
{\sc M. Thaddeus},
Stable pairs, linear systems and the Verlinde formula,
{\it Invent.\ Math.\ }117 (1994) 317--353.

\bibitem{thad2}
{\sc M. Thaddeus},
An introduction to the topology of the moduli space of stable bundles
on a Riemann surface,
{\it Geometry and physics (Aarhus, 1995)}, 71--99,
Lect. Notes in Pure and Appl. Math.\ 184, Dekker, 1997.

\bibitem{thad3}
{\sc M. Thaddeus},
On the Morse decomposition of the space of Higgs bundles on a curve,
in preparation.

\bibitem{zag}
{\sc D. Zagier},
On the cohomology of moduli spaces of rank two vector bundles over
curves,
{\it The moduli space of curves (Texel Island, 1994)}, 533--563,
Progr.\ Math., 129, Birkh\"auser, 1995.

}

\end{thebibliography}
\end{document}